\documentclass[12pt,reqno]{article}
mm \textheight=9.in

\usepackage{amssymb}
\usepackage{amsmath}
\usepackage{amsthm}
\usepackage{amssymb}
\usepackage{amsmath}
\usepackage{amsthm}
\usepackage{graphics}
\usepackage{color}
\newcommand{\tw}[3]{{$#1$}${\,\scriptscriptstyle {#2}}\atop\raise9pt\hbox{$\scriptstyle\tp$} ${$#3$}}
\newcommand{\st}[1]{\mbox{${\,\scriptscriptstyle {#1}}\atop\raise5.5pt\hbox{$*$}$}}
\newcommand{\btr}{\raise1.2pt\hbox{$\scriptstyle\blacktriangleright$}\hspace{2pt}}
\newcommand{\braid}[3]{{#1}$\lower4pt\hbox{$\tp\atop\raise4pt
            \hbox{$\scriptscriptstyle {#2} $}$}${#3}}

\newcommand{\id}{\mathrm{id}}

\newcommand{\Tc}{\mathcal{T}}

\newcommand{\A}{\mathcal{A}}
\newcommand{\B}{\mathcal{B}}

\newcommand{\Yg}{\mathrm{Y}}
\newcommand{\Ug}{\mathrm{U}}
\newcommand{\Xg}{\mathrm{X}}

\newcommand{\D}{\mathrm{D}}

\newcommand{\Mg}{\mathfrak{M}}
\newcommand{\Bg}{\mathfrak{B}}

\newcommand{\z}{\mathfrak{z}}

\renewcommand{\S}{\mathrm{S}}
\newcommand{\Ha}{\mathcal{H}}
\newcommand{\Ru}{\mathcal{R}}

\newcommand{\Sc}{\mathcal{S}}

\newcommand{\C}{\mathbb{C}}
\newcommand{\Z}{\mathbb{Z}}

\newcommand{\N}{\mathbb{N}}

\newcommand{\tp}{\otimes}

\newcommand{\zt}{\zeta}

\newcommand{\F}{\mathcal{F}}
\newcommand{\ve}{\varepsilon}
\newcommand{\gm}{\gamma}

\newcommand{\ot}{\otimes}

\newcommand{\tr}{\triangleright}
\newcommand{\tl}{\triangleleft}

\newcommand{\End}{\mathrm{End}}

\newcommand{\Res}{\mathrm{Res}}
\newcommand{\Hom}{\mathrm{Hom}}

\newcommand{\Tr}{\mathrm{Tr}}

\newcommand{\btl}{\mbox{\raise1.1pt\hbox{$\scriptstyle\blacktriangleleft$}}}

\newcommand{\g}{\mathfrak{g}}

\newcommand{\nn}{\nonumber}

\renewcommand{\l}{\mathfrak{l}}

\newcommand{\T}{\mathrm{T}}

\newcommand{\I}{\mathcal{I}}
\newcommand{\si}{\sigma}
\newcommand{\al}{\alpha}

\newcommand{\bt}{\beta}
\newcommand{\kp}{\kappa}
\newcommand{\be}{\begin{eqnarray}}
\newcommand{\ee}{\end{eqnarray}}

\newtheorem{thm}{Theorem}
\newtheorem{propn}[thm]{Proposition}
\newtheorem{lemma}[thm]{Lemma}
\newtheorem{corollary}[thm]{Corollary}

\theoremstyle{definition}
\newtheorem{remark}[thm]{Remark}
\newtheorem{definition}[thm]{Definition}

\newtheorem{parag}{}[section]

\begin{document}
\title{Reflection equation and twisted Yangians}
\author{A. I. Mudrov}
\date{}
\maketitle
\begin{center}
{
Department of Mathematics, University of York, YO10 5DD, UK,\\
 \texttt{aim501@york.ac.uk}}\\
 St.-Petersburg Department of Steklov Mathematical Institute,
\\
Fontanka 27, 191011 St.-Petersburg, Russia,\\
 \texttt{mudrov@pdmi.ras.ru}
\end{center}
\begin{abstract}
With any involutive anti-algebra and coalgebra automorphism
of a quasitriangular bialgebra we associate a reflection
equation algebra. A Hopf algebraic treatment
of the reflection equation of this type and its universal
solution is given. Applications to the twisted Yangians are considered.
\end{abstract}
\tableofcontents

\section{Introduction}
One purpose of the present work is to give a formal insight to some constructions
of the twisted Yangian theory from the Hopf algebraic point of view.
We believe that such an approach is useful for better understanding of
the  twisted Yangians and other algebras related to the so called reflection equation  (RE).
Another goal we pursue in this paper is to define and study a universal solution to the reflection
equation, with applications to the twisted Yangians in sight.

One of the key objects in the theory of quantum group is the  RTT equation, \cite{FRT},
associated a solution  $R\in \End(\C^{N})\tp \End(\C^{N})$
of the Yang-Baxter equation
$$
R_{12}R_{13}R_{23}=R_{23}R_{13}R_{12}
$$
in the algebra $\End(\C^{N})\tp \End(\C^{N})\tp \End(\C^{N})$.
Here we use the standard notation marking the tensor $\End(\C^N)$-factors with subscripts.
The RTT equation reads
\be
R_{12}T_1T_2&=&T_2T_1R_{12},
\label{rtt}
\ee
where $T=||T_{ij}||$ is $N\times N$-matrix with coefficients in some associative algebra. In particular,
$T_{ij}$ can be  the generators of
the RTT algebra, which is defined as a quotient of the free algebra $\C\langle T_{ij}\rangle$ by the relations
(\ref{rtt}).

An important counterpart of (\ref{rtt}) is the reflection equation,
which appeared in \cite{C} and found numerous applications in mathematical physics
and representation theory.
We distinguish the following two types of  RE:
\be
R_{21}K_1R_{12}K_2&=&K_2R_{21}K_1R_{12},
\label{RE1}\\
R_{12}S_1R'_{12}S_2&=&S_2(R')_{21}S_1R_{12},
\label{RE2}
\ee
see \cite{RS}, \cite{Mj1}, and \cite{KS}.
Here $K=||K_{ij}||$ and $S=||S_{ij}||$ are $N\times N$ matrices with coefficients in some associative algebras. The dashed operator
$R'$ is equal to $(\tau\tp \id)(R)$, where
$\tau$ is an involutive anti-automorphism of the algebra $\End(\C^N)$.
Similarly to the RTT case,  the RE algebra of relative type is defined as a quotient of the free algebra generated by
either $K_{ij}$ or $S_{ij}$ modulo the relations written as the corresponding RE.

Equations (\ref{RE1}) and (\ref{RE2})   have very much in common although
there are certain principal differences.
Consider the situation when $\C^N$ is a module over some quasitriangular Hopf algebra $\Ha$
and $R$ is the image of the universal R-matrix $\Ru\in \Ha\tp \Ha$.
Then the equation (\ref{RE1}) can be interpreted in terms of the Hopf algebra twist.
Namely, the algebra generated by the matrix elements $K_{ij}$ is a module over
the twisted tensor square \tw{\Ha}{\Ru}{\Ha}, which admits a Hopf algebra homomorphism
from the double $\D\Ha$. Important  in this picture is that $\Ha$ is endowed with
antipode. Then the corresponding RTT algebra, which is a natural $\Ha$-bimodule,
can be viewed as an algebra over  $\Ha^{op}\tp \Ha$, the trivial tensor product of the Hopf algebra
$\Ha$ and
its co-opposite.  The Hopf algebra \tw{\Ha}{\Ru}{\Ha} is a twist of
$\Ha^{op}\tp \Ha$.
Under this twist, the RE algebra is obtained from the RTT algebra associated with $R$, \cite{DM}.
This fact relates in many respects the theories of the equations (\ref{rtt}) and (\ref{RE1}),
e.g. allows to define a universal K-matrix, \cite{DKM}.

Equation (\ref{RE2}) appears in the theory of twisted Yangians and their applications, \cite{MNO}.
It turns out that this type  of RE can be given a similar treatment.
Now, however, we need not  assume   $\Ha$ to be a Hopf algebra.
The role of antipode can be given to an anti-algebra and coalgebra  automorphism $\tau\colon \Ha\to \Ha$,
which we assume to be involutive. Such an automorphism
makes the RTT algebra a module
over $\Ha\tp \Ha$. The latter bialgebra admits a twist to \tw{\Ha}{\Ru_\tau}{\Ha} with $\Ru_\tau:=(\tau\tp\tau)(\Ru^{-1})$.
This twist converts the RTT algebra to the RE algebra of the second type,
which thereby becomes a module algebra over \tw{\Ha}{\Ru_\tau}{\Ha}.
Note that $\Ru_\tau$ is a universal R-matrix of $\Ha$, as well as $\Ru$.
Thus the quasitriangular bialgebra $\Ha$ accounts for many properties of
the reflection equation (\ref{RE2}) relating them to properties of the RTT algebra.
In particular, it makes  rather natural definition of  the universal solution of (\ref{RE2})
or universal S-matrix.

For a  fixed associative algebra $\A$, we are looking for an element from $\Ha\tp \A$ satisfying (\ref{RE2}) in
the tensor product $\Ha\tp \Ha\tp \A$,
with $R$ replaced by  $\Ru$ and the subscripts referring to the tensor $\Ha$-factors.
Given such an element, any representation of $\Ha$ yields a matrix solution with coefficients in $\A$.
For instance, we can take the field of scalars for the role of $\A$;
then the solution in question will be an element from $\Ha$. In any representation of $\Ha$,
it produces a numerical  RE matrix.

A universal solution of (\ref{RE2}) can be interpreted as a collection of intertwiners between certain modules
over the RE algebra. Numerical RE matrices arising from characters of the RE algebra appear in a recent work
\cite{KN}
as the images of the Zhelobenko cocycles \cite{Zh}, under a functor between the module categories
of the orthogonal (symplectic) Lie algebra and of twisted Yangians.
In our theory those intertwiners become a part of a more general structure,
namely, the universal S-matrix.

Specifically, we define the universal S-matrix as an element $\Sc\in \Ha\tp \A$ satisfying the equation
\be
\label{characteristic0}
(\Delta\tp \id)(\Sc)=\Sc_{13}\Ru'_{12}\Sc_{23}, \quad\mbox{with} \quad \Ru':=(\tau\tp\id)(\Ru),
\ee
in $\Ha\tp \Ha\tp \A$.
The intertwining property of $\Ru$
guarantees that $\Sc$ solves (\ref{RE2}), hence $\Sc$ is indeed a universal solution of (\ref{RE2}).
The element $\Sc$ depends on $\A$ but single $\A$ may give different S-matrices.
Roughly speaking, they are parameterized by homomorphisms of the RE algebra to $\A$.
If, $\A$ can be constructed on the vector space $\Ha^*$; then it admits
a non-degenerate pairing with $\Ha$.
In this case $\S$ is the canonical element of the paring
or the identity operator upon the identification $\Ha\tp \Ha^*\simeq \End(\Ha)$.

Here is  one more explanation why we put (\ref{characteristic0}) to be the central equation of our theory.
We can regard $\Sc$ as a collection of matrix solutions of  (\ref{RE2}) labeled by representations of $\Ha$.
This collection forms a category, whose homsets consist of  $\Ha$-equivariant map "commuting" with S-matrices.
The equation  (\ref{characteristic0})  facilitates
a fusion procedure or, in other words, a monoidal structure
on that category.

It turns out that the above picture is too simple to incorporate the twisted Yangian $\Yg(\g_N)$,
which is  manifestly defined through an equation of the type (\ref{RE2}), see \cite{MNO}.
Here $\g_N\subset \g\l_N$ is the Lie algebra preserving either symmetric or skew symmetric non-degenerate bilinear form.
The twisted Yangian is a comodule over the Yangian $\Yg(\g\l_N)$ and hence a module over its dual $\Yg^*(\g\l_N)$.
However $\Yg^*(\g\l_N)$ is not a quasitriangular bialgebra.
By this reason,  one cannot obtain  $\Yg(\g_N)$ from $\Yg(\g\l_N)$, which is a sort of  RTT
algebra, directly applying the technique of twist. The question is
what quasitriangular bialgebra should be taken as $\Ha$.

At a first glance, that could be the double $\D\>\Yg^*(\g\l_N)$, which
is a quasitriangular Hopf algebra and hence a natural candidate for the role of $\Ha$.
The Hopf dual to $\D\>\Yg^*(\g\l_N)$ has a projection to $\Yg(\g\l_N)$ making
the corresponding RE algebra a comodule over  $\Yg(\g\l_N)$.
One  might expect that  $\Yg(\g_N)$ could be obtained by a homomorphism
from the RE dual of  $\D\>\Yg^*(\g\l_N)$. For example, that  homomorphism might be implemented by means of an RE character,
as $\Yg(\g_N)$ is a subalgebra in  $\Yg(\g\l_N)$.
However this is not the right way to go. Quite surprisingly, the
RE dual to the  double of an involutive Hopf algebra admits no characters at all.
At the same time, it is known that  any symmetric or skew-symmetric numerical
$N\times N$ matrix defines a one dimensional representation of $\Yg(\g_N)$.

It turns out, in what concerns the twisted Yangians,  that the role of $\Ha$
should be given to another (pseudo) quasitriangular bialgebra, which nevertheless
admits a non-trivial homomorphism from the double. Thus all the modules entering  $\Sc$
are still   $\D\>\Yg^*(\g\l_N)$- and hence  $\Yg(\g\l_N)$-modules. They include, in particular,
the modules associated with skew Young diagrams, \cite{NT}.

The setup of the paper is as follows. The next section contains basic facts from quantum groups.
Section \ref{Sec:gen} develops the general theory of the reflection equation (\ref{RE2}).
It contains the definition of RE twist, RE dual and its symmetries, the universal S-matrix and its basic properties.
Section \ref{Sec:QD} is a detailed consideration of the special case when $\Ha$ is a quantum double.
The main result of this section is factorization of the $S$-matrix of the double.
In Section \ref{Sec:FRTform} we introduce reduced RE algebras,  whose symmetries do not form a
quasitriangular bialgebra. We do it using the dual version of the Faddeev-Reshetikhin-Takhtajan formalism.
In Section \ref{Sec:Pseudo} we define the universal S-matrix for
pseudo-quasitriangular bialgebras, whose R-matrix is not invertible.
Section \ref{Sec:TY} deals with applications of the general theory to the twisted Yangians.

\vspace{0.2in}
\noindent
{\bf \large Acknowledgements}
This research is supported by the EPSRC grant C511166 and partly
supported by the RFBR grant 06-01-00451.
We are grateful to Maxim Nazarov for formulating the problem and for valuable discussions.
Our special thanks are to the referee for his meticulous reading of the manuscript.
His numerous remarks and suggestions have served to considerable improvement of the text.
\section{Some basics on quasitriangular bialgebras}
\label{sGTRE}
\begin{parag}{\bf Quasitriangular bialgebras.}
\label{QHA}
The reader is assumed to be familiar with the general theory of Hopf algebras, quasitriangular structure and
twist. For a guide in this field, we refer to any textbook in quantum groups, e.g.
\cite{ChP,Mj2}. Some topics concerning the twisted tensor product of Hopf algebras are discussed in \cite{DM}.

For the sake of simplicity we assume in Section \ref{sGTRE}--\ref{Sec:QD} that the bialgebra $\Ha$ is {\em finite dimensional}.
The generalization to infinite dimensional case requires certain completion of tensor products. We consider the
necessary modifications for a particular infinite dimensional bialgebra in Section \ref{Sec:FRTform}.

Consider a quasitriangular bialgebra $\Ha$ over the complex field with the coproduct $\Delta$, counit $\ve$,
and universal R-matrix $\Ru$, \cite{D2}.
Throughout the paper we adopt the following convention concerning the presentation of tensors. We use the Sweedler symbolic notation
$h^{(1)}\tp h^{(2)}:=\Delta(h)$ for the coproduct of an element $h\in \Ha$. For more general tensors we mark their components with subscripts
and assume
suppressed summation, like $\Ru=\Ru_1\tp \Ru_2\in \Ha\tp\Ha$.
To distinguish the Sweedler components of different copies of the same tensor, we use primes,
e.g. $\Ru\Ru\Ru_{21}=\Ru_1\Ru_{1'}\Ru_{2''}\tp \Ru_2\Ru_{2'}\Ru_{1''}\in\Ha\tp\Ha$.

Recall that a quasitriangular structure  $\Ru$  on $\Ha$ (universal R-matrix) is an invertible element from $\Ha\tp \Ha$
 that satisfies the identities
\be
\label{hex}
(\Delta\tp \id)(\Ru)=\Ru_{13}\Ru_{23},\quad (\id \tp \Delta)(\Ru)=\Ru_{13}\Ru_{12},
\ee
\be
\Ru\Delta(h)=(\si\circ\Delta)(h)\Ru
\ee
for all $h\in \Ha$. Here $\si$ is the flip of the tensor factors, $\si\colon x\tp y\mapsto y\tp x$.

If $\Ha$ is a Hopf algebra, i.e. a bialgebra with antipode $\gm$, then
the R-matrix possesses the following symmetries:
\be
\label{GmR}
(\gm\tp \gm)(\Ru)=\Ru,\quad (\gm\tp \id)(\Ru)= \Ru^{-1},\quad (\id\tp \gm)(\Ru^{-1})=\Ru.
\ee
We consider only Hopf algebras with  invertible antipode.

The element $\Ru_{21}^{-1}$ is an alternative quasitriangular structure on
$\Ha$, which differs from $\Ru$ in general. In what follows, we use the notation
$$
\Ru^+:=\Ru,\quad\Ru^-:=\Ru^{-1}_{21}.
$$
It is also convenient to consider $\Ru^\pm$ as linear maps from the dual bialgebra
$\Ha^*$ to $\Ha$ acting by $a\mapsto \langle a,\Ru^\pm_1\rangle\Ru^\pm_2$ (pairing with the first tensor component of $\Ru$). These maps establish bialgebra
homomorphisms from $\Ha^{*op}$ to $\Ha$. Here and further on
the superscript $op$ designates the bialgebra with the opposite comultiplication. The subscript
$op$ stands for the opposite multiplication.

\end{parag}
\begin{parag}{\bf Twisted tensor square.}
\label{Ssec:TTS}
Let $\F\in \Ha\tp \Ha$ be an invertible element  obeying the conditions
\be
&(\Delta\tp \id)(\F)\F_{12}=(\id\tp \Delta)(\F)\F_{23},
\label{twist}\\
&(\ve\tp \id)(\F)=1\tp 1=(\id\tp \ve)(\F).
\nn
\ee
Here the first equation is understood in $\Ha^{\tp 3}$ while the second in $\Ha^{\tp 2}$, upon
the embedding $1\colon\C\to \Ha$. The element $\F$ is called twisting cocycle.

Recall that  twist of a bialgebra $\Ha$ is a new bialgebra $\tilde\Ha$ with
the same multiplication and counit, but the new comultiplication
$$
\tilde \Delta(h)=\F^{-1}\Delta(h)\F.
$$
If $\Ha$ is a (quasitriangular) Hopf algebra with the universal R-matrix $\Ru$, then $\tilde\Ha$ is a (quasitriangular) Hopf algebra as well,
\cite{Dr3}.

Given a twisting cocycle $\F$, we denote by $\F^{(k)}\in \Ha^{\tp k}$, $k=0,1,\ldots$, the intertwining operator between the twisted and non-twisted
$k$-folded coproducts, $\F^{(k)}\tilde\Delta^k(h)=\Delta^k(h)\F^{(k)}$, for all $h\in \Ha$.
This operator can be constructed inductively using the factorization $\F^{(n)}=(\Delta^k\tp \Delta^m)(\F)\>(\F^{(k)}\tp \F^{(m)})$,
where $n=k+m$, $\F^{(0)}=1_\C$, $\F^{(1)}=1_\Ha$, and $\F^{(2)}=\F$. The right-hand side is independent on the partition $n=k+m$, as follows from (\ref{twist}).

A particular case of twist called twisted tensor product will be of special interest for us.
We recall the corresponding definition for the  twisted tensor square of a quasitriangular bialgebra.
Consider the universal R-matrix of $\Ha$ as an element $\Ru_{23}$
from $(\Ha\tp \Ha)^{\tp 2}$. Then it is a twisting cocycle of the bialgebra
$\Ha\tp \Ha$ (the latter is naturally equipped with  algebra and coalgebra structures).
We define \tw{\Ha}{\Ru}{\Ha} to be the twisted bialgebra coinciding with $\Ha\tp \Ha$ as
an associative algebra and endowed with the comultiplication
$$
\Delta(x\tp y)=\Ru^{-1}_{23}\Bigl(\bigl(x^{(1)}\tp y^{(1)}\bigr)\tp\bigl( x^{(2)}\tp y^{(2)}\bigr)\Bigr) \Ru_{23}.
$$
In the right-hand side, the Sweedler notation  stands for the coproduct in $\Ha$.
Note that one can define the twisted tensor square with $\Ru$ replaced by any invertible element from $\Ha\tp \Ha$ obeying (\ref{hex}).

Since $\Ha\tp \Ha$ can be equipped with a quasitriangular structure
(generally speaking, in various ways), that also holds for \tw{\Ha}{\Ru}{\Ha}.
It is known that the algebra maps
$$
\Ha\stackrel{\Delta}{\longrightarrow} \Ha\tp \Ha,\quad \Ha^{*op}\stackrel{\Delta}{\longrightarrow}\Ha^{*op}\tp \Ha^{*op}
\stackrel{\Ru^+\tp \Ru^-}{\longrightarrow}  \Ha\tp \Ha
$$
induce bialgebra homomorphisms
\be
\Ha\to \mbox{\tw{\Ha}{\Ru}{\Ha}},\quad
\Ha^{*op}\to \mbox{\tw{\Ha}{\Ru}{\Ha}}.
\label{morphisms}
\ee
If $\Ha$ is a Hopf algebra, then one can define the double $\D\Ha$,
which is also a quasitriangular Hopf algebra containing $\Ha$ and $\Ha^{*op}$ as
Hopf subalgebras. Then the maps (\ref{morphisms}) extend to a Hopf algebra homomorphism from
$\D\Ha$ to \tw{\Ha}{\Ru}{\Ha}.

We  denote by $\braid{\Ha^*}{\Ru}{\Ha^*}$ the dual algebra to \tw{\Ha}{\Ru}{\Ha}.
It contains two copies of $\Ha^*$ both embedded as sub-bialgebras.
There are bialgebra homomorphisms
\be
\mbox{$\braid{\Ha^*}{\Ru}{\Ha^*}$}\to \Ha_{op},
\quad
\mbox{$\braid{\Ha^*}{\Ru}{\Ha^*}$}\to \Ha^*
\label{dualmaps}
\ee
obtained by dualization of the arrows (\ref{morphisms}).
They act by
\be
x\tp y\mapsto \Ru^{-1}(y)\Ru_{21}(x),
\quad
x\tp y\mapsto xy,
\label{dualmaps_exp}
\ee
respectively. Note that  the right-hand side in the left formula is expressed in terms of the multiplication
in $\Ha$ rather than $\Ha_{op}$.

\end{parag}
\begin{parag}{\bf Involutive automorphisms of bialgebras.}
In this subsection we consider endomorphisms of the bialgebra
$\Ha$ preserving the multiplication and comultiplication modulo permutation of the factors.
The term automorphism is reserved for  algebra and coalgebra invertible maps $\Ha\to \Ha$.

Denote by $\Z_2=\{-,+\}$ the group of two elements assuming $+$ to be the identity.
Let $\Ha$ be a bialgebra and let $\tau\colon\Ha\to \Ha$ be a linear endomorphism.
Suppose that  $\tau$ preserves the multiplication and comultiplication modulo permutation
of the factors. We say that $\tau$ has {\em signature} $\mu_1\choose\mu_2$, where $\mu_1,\mu_2\in \Z_2$,
setting
\be
\mu_1&=&
\left\{
\begin{array}{cll}
+,&\tau \mbox { is a coalgebra map}\\
-,&\tau \mbox { is an anti-coalgebra map}
\end{array}
\right.
,\quad
\mu_2=
\left\{
\begin{array}{cll}
+,&\tau \mbox { is an algebra map}\\
-,&\tau \mbox { is an anti-algebra map}
\end{array}
\right.
\nn
\ee
We call $\tau$ {\em skew} if its signature $\mu_1\choose \mu_2$ satisfies the condition
$\mu_1\mu_2=-$.

If two endomorphisms have definite signature, so does its composition, and the
corresponding signs are multiplied. If an endomorphism has signature
 $\mu_1\choose \mu_2$  then the dual conjugate endomorphism has signature
  $\mu_2\choose \mu_1$.
\begin{propn}
\label{sign}
Let $\tau$ be a $\mu_1\choose\mu_2$-automorphism of a Hopf algebra $\Ha$
with the counit $\ve$ and invertible antipode $\gm$.
Then a) $\ve\circ \tau=\ve$, b)
 $\tau\circ \gm\circ \tau^{-1}=\gm^{\mu_1\mu_2}$, and
c) if $\tau$ is a skew  involution, so is $\tau\circ \gm^{k}$ for all $k\in \Z$.
\end{propn}
\begin{proof}
Uniqueness of the counit implies a).
Obviously c) follows from b). To check b), observe that the map
and $\tau\circ \gm\circ \tau^{-1}$
satisfies the axioms for the antipode in $\Ha^{op}$. But so does $\gm^{-1}$, hence
b) follows from uniqueness of the antipode.
\end{proof}
\end{parag}
\section{General theory of reflection equation}
\label{Sec:gen}
\begin{parag}{\bf Reflection equation twist.}
Let $\Ha$ be a quasitriangular bialgebra with the universal R-matrix $\Ru$.
Regard the inverse $\Ru^{-1}$ as a universal R-matrix for the opposite bialgebra $\Ha_{op}$ (with the opposite multiplication).
Then the tensor product $\Ha_{op}\tp \Ha$ is a quasitriangular bialgebra with the universal R-matrix
$\Ru^{-1}_{13}\Ru_{24}$. Modules over $\Ha_{op}\tp \Ha$ are natural bimodules over $\Ha$ and {\em vice versa}.
Suppose that $\A$ is an $\Ha$-bimodule algebra and suppose $\A$ is commutative in the category of
$\Ha$-bimodules. Commutativity means the equality
\be
(a\tl\Ru_1)(b\tl\Ru_2)=(\Ru_2\tr b)(\Ru_1\tr a)
\label{RTT}
\ee
for all $a,b\in \A$, where $\tr$ and $\tl$ are respectively the left and right $\Ha$-actions.
They are obtained by duality  from the  right and left  regular actions of $\Ha$ on itself and hence commute with each other.
Recall that the subscripts mark the tensor components, as adopted throughout the paper.
One may think of $\A$ as the RTT algebra or its suitable quotients such as
the dual algebra $\Ha^*$; then the $\Ha$-actions $\tr$ and $\tl$ are dual to the right and left regular actions.

Suppose that  $\tau\colon \Ha\to \Ha$ is involutive $+\choose-$-automorphism.
The elements $\Ru^\pm_{\tau}:=(\tau\tp \tau)(\Ru^\pm)^{-1}$ are also universal R-matrices
of $\Ha$ (which in general do not coincide with $\Ru^\pm$.
In the presence of $\tau$, the $\Ha$-bimodule $\A$ becomes an $\Ha\tp \Ha$-module
with respect to the action
\be
(x\tp y)\vdash a= y\tr a \tl\tau (x).
\ee
 It is
a commutative algebra in the category of modules over $\Ha\tp \Ha$,
the latter being equipped with the quasitriangular structure $(\Ru_\tau)_{13}\Ru_{24}$.
This fact is expressed by the identity (\ref{RTT}).

Consider the bialgebra algebra twist from $\Ha\tp \Ha$ to \tw{\Ha}{\Ru_\tau}{\Ha}.
We call this transition  the RE twist, for the reason to be clear later on.
Any  twist induces a transformation of module algebras; they are endowed with
new multiplications that are consistent with the new coproduct in the twisted bialgebra.
As applied to the $\Ha\tp \Ha$-module algebra $\A$, the RE twist produces a new
algebra, $\tilde \A$,
over \tw{\Ha}{\Ru_\tau}{\Ha}. The multiplication $\tilde m$ in $\tilde \A$ is related with the multiplication
$m$ in $\A$ by the formulas
\be
\tilde m(a\tp b)= m\bigl(\tau(\Ru^{-1}_1)\tr a\tp b\tl\Ru^{-1}_2\bigr),
\quad
m(a\tp b)= \tilde m\bigl(\tau(\Ru_1)\tr a\tp b\tl\Ru_2\bigr).
\label{prod}
\ee
As usually, we use the subscripts to denote the tensor components with implicit summation:
$\Ru^{-1}=\Ru^{-1}_1\tp \Ru^{-1}_2$.
Expressing the multiplication $m$  through $\tilde m$ in the formula (\ref{RTT}),
we come to the corresponding relation in $\tilde \A$:
\be
\bigl(\tau(\Ru_{1'})\tr a\tl\Ru_1\bigr)(b\tl\Ru_2 \Ru_{2'})=\bigl(\tau(\Ru_{1'})\Ru_2\tr b\bigr)(\Ru_1\tr a\tl \Ru_{2'}).
\label{REu'}
\ee
Since the algebra $\A$ is commutative in the category of $\Ha\tp \Ha$-modules,
the algebra $\tilde \A$ is commutative in the category of
\tw{\Ha}{\Ru_\tau}{\Ha}-modules, which fact is expressed by the equation (\ref{REu'}).
\begin{remark}
Note that $\A$ may have other identities besides (\ref{RTT}). They can also be expressed in
terms of the new multiplication, producing identities in $\tilde \A$. The equation
(\ref{REu'}) is a form of the "reflection equation",
as (\ref{RTT}) is a form of the RTT relation.
\end{remark}
The most important example of $\A$ is the dual algebra $\Ha^*$.
In this case, the algebra $\tilde \A$ will be denoted by $\Ha^\circledast$.
For infinite dimensional bialgebras considered in Section \ref{Sec:FRTform},
$\A$ is a restricted  dual algebra to $\Ha$, that is,
an $\Ha$-bimodule subalgebra $\Ha^*$ in the convolution algebra of all linear functions
on $\Ha$, whose natural paring with $\Ha$ is nondegenerate.

\end{parag}

\begin{parag}{\bf Standard form of RE.}
Consider an associative algebra   $\Mg$ and an epimorphism
$\rho\colon \Ha \to \Mg$. Suppose that the dual space $\Mg^*$ is embedded in the  dual algebra $\Ha^*$
and generates it.
Let $R$ denote the image of the universal R-matrix in  $\Mg\tp \Mg$.
Suppose there exists an involutive anti-automorphism of $\Mg$ intertwined with $\tau$
by $\rho$. We will denote it with the same letter $\tau$; so by definition $\tau\circ\rho=\rho\circ\tau$.

Let $\{e_i\}\subset \Mg$ be a basis and let  $\{e^i\}\subset \Mg^*$ be its dual.
Denote by $S$ the element $\sum_i e_i\tp e^i$.
With $a,b$ set to be $e_i,e^j$ the relation (\ref{REu'}) turns into the matrix reflection equation
\be
R_{12} S_1 R'_{12}S_2&=&S_2R''_{12}S_1 R_{12} ,
\label{RE}
\ee
where  $R'=\tau(R_1)\tp R_2$ and $R''=R_2\tp \tau(R_1)$,
so that $R''=(R')_{21}$.
We shall use the similar notation for the universal R-matrices.
It may happen  that $(\tau\tp \tau)(R)=R_{21}$ and therefore $R''=R'$, as is the case with
the Yang matrix, see the corresponding section.
\end{parag}
\begin{parag}{\bf Universal RE matrix.}
\label{universal}
We introduce the universal matrices $\Tc\in \Ha\tp \Ha^*$ and  $\Sc\in \Ha\tp \Ha^\circledast$ by setting
$\Tc:=\sum_\mu e_\mu\tp e^\mu=:\Sc$, where $\{e_\mu\}$ is a basis in $\Ha$ and $\{e^\mu\}$ is the dual basis
in $\Ha^*$ (recall that $\Ha^*$ coincides with $\Ha^\circledast$ as a vector space).

The T-matrix fulfills the identity
\be
(\Delta\tp \id)(\Tc)=\Tc_1\Tc_2,
\label{charateristicT}
\ee
in $\Ha\tp\Ha\tp  \Ha^*$, which implies the "universal" RTT relation
\be
\Ru_{12}\Tc_1\Tc_2=\Tc_2\Tc_1\Ru_{12}.
\label{characteristicT}
\ee
Here the subscripts mark the $\Ha$-tensor factors.
Expressing the multiplication in $\Ha^*$ in terms of the multiplication in $\Ha^\circledast$  by  the formula (\ref{prod})
we obtain the equation
\be
(\Delta\tp \id)(\Sc)=\Sc_1\Ru'_{12}\Sc_2
\label{characteristic}
\ee
held in $\Ha\tp \Ha\tp \Ha^\circledast$. As above, the subscripts indicate the $\Ha$-tensor
factors.
From (\ref{characteristicT}) we conclude that $\Sc$ satisfies the equation
\be
\label{uRE}
\Ru\Sc_1\Ru'\Sc_2=\Sc_2\Ru''\Sc_1 \Ru.
\ee
This fact justifies the following definition.
\begin{definition}
Let $\Ha$ be a  quasitriangular bialgebra with involution.
A universal S-matrix over $\Ha$ with coefficients in an associative algebra $\A$ is
a non-zero element $\Sc\in \Ha\tp \A$ obeying the identity (\ref{characteristic}) in $\Ha\tp \Ha\tp \A$.
\end{definition}
We remark that "universality" in this definition is relevant to the $\Ha$-component of $\Sc$,
because there can be different S-matrices with the same $\A$.
The case of $\A=\Ha^\circledast$ with $\Sc$
the canonical element is distinguished, as any other S-matrix can  be obtained from it, see Proposition \ref{hom} below.
Hence the canonical element of $\Ha\tp \Ha^\circledast$ can be called the universal S-matrix
without referring to coefficients.
\begin{propn}
\label{hom}
A universal S-matrix $\Sc_\A$ with coefficients in $\A$ defines a ring homomorphism
$\phi\colon \Ha^\circledast\to \A$ such that $\Sc_\A=(\id\tp \phi)(\Sc)$, where $\Sc$ is the canonical element
of $\Ha\tp \Ha^\circledast$.
\end{propn}
\begin{proof}
Consider $\Sc_\A$ as a linear mapping  $\phi\colon \Ha^\circledast\to \A$ through
the pairing of $\eta\in \Ha^\circledast$ with the $\Ha$-component of $\Sc_\A$.
Here we use the identification $\Ha^\circledast=\Ha^*$. The
identity (\ref{characteristic}) implies for all $\eta,\xi\in \Ha^\circledast$
the equation  $\phi(\eta \xi)=\phi(\Ru'_1\tr\eta)\phi(\xi\tl\Ru'_2)$,
where the product $\eta \xi$ is taken in $\Ha^*$.
Using the right formula in (\ref{prod}) we express the multiplication in $\Ha^*$
through the multiplication in $\Ha^\circledast$. This shows that $\phi$ is multiplicative (here we also use invertibility of  $\Ru$).
Applying $\phi$ to the $\Ha^\circledast$-component of the canonical element, we obtain $\Sc_\A$.
\end{proof}

\end{parag}
\begin{parag}{\bf Elementary properties of universal S-matrix.}
We assume that the algebra $\A$ is integral over its center $\z(\A)$, that is,
for all $z\in \z(\A)$ and $a\in \A$ the equality $za=0$ implies either
$z=0$ or  $a=0$. Let $\Sc\in \Ha\tp \A$ be a universal S-matrix, where $\Ha$ is a bialgebra with involution.
We consider $\Sc$ as a map $\Ha^\circledast\to \A$ implemented through the paring between $\Ha$ and $\Ha^\circledast$
and suppose that  $\A$ is the image of $\Ha^\circledast$ (cf. Proposition \ref{hom}).
\begin{lemma}
\label{S-counit}
Under the above assumptions, $(\ve\tp \id)(\Sc)=1$.
\end{lemma}
\begin{proof}
Applying the counit to the left (right) $\Ha$-components of (\ref{characteristic}) we
get
$$(\ve\tp \id)(\Sc)\Sc=\Sc=\Sc(\ve\tp \id)(\Sc).$$
Hence $(\ve\tp \id)(\Sc)$ belongs to $\z(\A)$ and is equal to $1$, since $\Sc\not=0$ by definition.
\end{proof}
The following assertion is also formulated under the above assumptions.
\begin{propn}
\label{char_inverse}
Suppose that $\Ha$ is a Hopf algebra with the antipode $\gm$.
Then the matrix $\Sc$ in invertible in $\Ha\tp \A$ with the  inverse
$$
\bigl((\gm\tau)(\Ru_1)\tp  1\bigr)(\gm\tp \id)(\Sc)\bigl(\Ru_2\tp 1\bigr)=\Sc^{-1}=\bigl(\tau(\Ru_1)\tp 1\bigr)(\gm\tp \id)(\Sc)\bigl(\gm(\Ru_2)\tp 1\bigr).
$$
\end{propn}
\begin{proof}
Applying the antipode to the left (right) $\Ha$-component of (\ref{characteristic}) and then multiplying the components, we get
$$
\bigl((\gm\tau)(\Ru_1)\tp  1\bigr)(\gm\tp \id)(\Sc)\bigl(\Ru_2\tp 1\bigr)\Sc=(\ve\tp \id)(\Sc)=\Sc\bigl(\tau(\Ru_1)\tp 1\bigr)(\gm\tp \id)(\Sc)\bigl(\gm(\Ru_2)\tp 1\bigr).
$$
This proves the statement, in view of Lemma \ref{S-counit}.
\end{proof}

We define
the element $\Sc^{op}:=(\bar\Ru_2\tp 1)\Sc\bigl(\tau(\bar\Ru_1)\tp 1\bigr)\in \Ha\tp \A$, where $\bar \Ru:=\Ru^{-1}$.
Now we do not require the antipode be present in $\Ha$.
\begin{propn}
The element $\Sc^{op}$ is a universal S-matrix
for the co-opposite bialgebra $\Ha^{op}$ with the R-matrix $\bar\Ru$.
\end{propn}
\begin{proof}
Equations (\ref{hex}) and (\ref{characteristic}) imply
$$
(\Delta\tp \id)(\Sc^{op})=\bigl(\bar\Ru_{2'}\bar\Ru_{2}\tp \bar\Ru_{2'''}\bar\Ru_{2''}\bigr)\bigl(\Sc_1\Ru'_{12}\Sc_2\bigr)
\bigl(\tau(\bar\Ru_{1''})\tau(\bar\Ru_{1})\tp\tau(\bar\Ru_{1'''}) \tau(\bar\Ru_{1'})\bigr),
$$
where
the primes distinguish different copies of R-matrices.
Now notice that $\Ru'$ in the middle factor cancels $\Ru$ whose components are labeled with the double prime. What is left gives
the identity
$$
(\Delta_{op}\tp \id)(\Sc^{op})=\Sc^{op}_1 \bar\Ru_{12}'\Sc^{op}_2,
$$
as required.
\end{proof}
\end{parag}
\begin{parag}{\bf Conjugate reflection equation.}
In some applications one has to work with different RE simultaneously.
Introduce the notation $(\Ru^{-1})':=(\tau\tp \id)(\Ru^{-1})$, $(\Ru^{-1})'':=(\id\tp \tau)(\Ru^{-1}_{21})$ and
define the matrix $\Sc^\bullet\in \Ha\tp\A$ to be a solution of the equation
\be
(\Delta_{op}\tp \id)(\Sc^\bullet)=\Sc^\bullet_1 (\Ru^{-1})''_{12} \Sc^\bullet_2.
\label{char_d}
\ee
This element satisfies the {\em conjugate} reflection equation
\be
\Ru^{-1}_{12} \Sc^\bullet_1 (\Ru^{-1})''_{12}\Sc^\bullet_2 &=& \Sc^\bullet_2(\Ru^{-1})'_{12}\Sc^\bullet_1 \Ru^{-1}_{12}.
\label{REdual}
\ee
The symmetry of this equation is the bialgebra \tw{\Ha^{op}}{\Ru_{21}}{\Ha^{op}}
(note that $\Ru_{21}$ is the universal $R$-matrix for the co-opposite Hopf algebra  $\Ha^{op}$).
The corresponding RE algebra is obtained from the RTT algebra via the twist
of $\Ha\tp \Ha$ by the cocycle $\Ru_{32}$.

The relation between the RE and its conjugate for $\Ha$ being a Hopf algebra is described by the following proposition.
\begin{propn}
Suppose that $\Ha$ is a Hopf algebra with involution $\tau$.  Let $\Sc^\bullet\in \Ha\tp\A$ be a universal solution
of the conjugate RE.
Then $\bar \Sc^\bullet:=(\Sc^\bullet)^{-1}$ is a universal S-matrix for the Hopf algebra $\Ha$ with
involution ${\tau\gm^2}$.
\end{propn}
\begin{proof}
Indeed, permute the tensor factors in the equality (\ref{char_d}) and take the inverse.
That results in the equality
$
(\Delta\tp \id)(\bar\Sc^\bullet)= \bar\Sc^\bullet_1(\bar\Ru')^{-1}\bar\Sc^\bullet_2,
$
where $\bar \Ru:=\Ru^{-1}$.
We would prove the statement if we show that
$
(\bar\Ru')^{-1}=(\tau\circ\gm^2\tp \id)(\Ru).
$
But the left-hand side of the latter equality is equal to
$
\bigl((\tau\circ\gm\tp \id)(\Ru)\bigr)^{-1}=(\tau\circ\gm\tp \id)(\Ru^{-1}),
$
because $\tau\circ\gm$ is an algebra automorphism.
This gives the right-hand side, due to   (\ref{GmR}).
\end{proof}
\end{parag}
\begin{parag}{\bf Comodule structures on the RE algebra.}
\label{Subsec:Sintertwiner}
In the present section we fix the domain of coefficients to be $\A=\Ha^\circledast$.
The algebra $\Ha^\circledast$ is a right comodule over
the dual bialgebra to
\tw{\Ha}{\Ru_\tau}{\Ha},
which we denote by  $\braid{\Ha^*}{\Ru_\tau}{\Ha^*}$. The coaction is given by the assignment
$$
\Ha^\circledast\ni a\mapsto a^{(2)}\tp \tau^*(a^{(1)})\tp a^{(3)}\in \Ha^\circledast\tp \mbox{$\braid{\Ha^*}{\Ru_\tau}{\Ha^*}$},
$$
and it is an algebra homomorphism.
Using  homomorphisms from $\braid{\Ha^*}{\Ru_\tau}{\Ha^*}$ to other bialgebras,
e.g. as in (\ref{dualmaps}), one can obtain various
comodule algebra structures over those bialgebras.

Let us consider in some detail the case of $\Ha^*$, using the left map from (\ref{dualmaps}),
where $\Ru$ is replaced by $\Ru_\tau$.
This way we obtain a right $\Ha^*$-coaction
$$
\Ha^\circledast\ni a\mapsto \delta^r_{\Ha^*}(a):=a^{(2)}\tp \tau^*(a^{(1)})a^{(3)}\in \Ha^\circledast\tp \Ha^*.
$$

Taking into account that the dual conjugate involution
$\tau^*\colon \Ha^*\to \Ha^*$ has signature $-\choose+$, we
have the following right and left $\Ha^*$-coalgebra structures on $\Ha^\circledast$:
\be
(\id\tp \delta^r_{\Ha^*})(\Sc)=\tau^*(\Tc_{13})\Sc_{12}\Tc_{13}
,\quad
(\id\tp \delta^l_{\Ha^*})(\Sc)=\Tc_{12}\Sc_{13}\tau^*(\Tc_{12}),
\label{dual_coaction}
\ee
written in terms of the S-matrix.
The bialgebra homomorphism $\Ru^+\colon \Ha^*\to \Ha^{op}$ converts these coactions into the right and left $\Ha^{op}$-coactions
\be
(\id\tp \delta^r_{\Ha^{op}})(\Sc)= \Ru_{13}' \Sc_{12} \Ru_{13}
,\quad
(\id\tp \delta^l_{\Ha^{op}})(\Sc)=\Ru_{12}\Sc_{13} \Ru_{12}'.
\ee

Further, consider the mappings
$$
\Phi^+\colon\eta\mapsto \langle\eta,\tau(\Ru_{1})\rangle\tau(\Ru_{2}),
\quad
\Phi^-\colon\eta\mapsto \langle\eta,\Ru_{2}\rangle\Ru_{1}
$$
defined on $\Ha^*_{op}$ with values in $\Ha$. It is easy to check that $\Phi^\pm$ are bialgebra homomorphisms.
\begin{propn}
\label{H*DH1}
The mapping $\Phi\colon \Ha^*_{op}\ni\eta\mapsto \Phi^+\bigl(\eta^{(1)}\bigr)\tp \Phi^-\bigl(\eta^{(2)}\bigr)\in \Ha\tp \Ha$
implements a bialgebra homomorphism $\Ha^{*}_{op}\to \mbox{\tw{\Ha}{\Ru_\tau}{\Ha}}$.
\end{propn}
\begin{proof}
By construction, $\Phi$ is an algebra map. Let us check that it preserves comultiplication.
The verification reduces to the formula
$$
\Phi^-(h^{(1)})\tp \Phi^+(h^{(2)})=(\Ru_\tau)^{-1}\Bigl(\Phi^-(h^{(2)})\tp \Phi^+(h^{(1)})\Bigr)\Ru_\tau
$$
for all $h\in \Ha$. This is equivalent to the identity
$$
(\Ru_\tau)_{13}\Ru^-_{23}(\Ru^{-1}_\tau)_{12}=(\Ru^{-1}_\tau)_{12}\Ru^-_{23}(\Ru_\tau)_{13},
$$
which is clearly true.
\end{proof}

The bialgebra mapping from Proposition \ref{H*DH1} induces a left action
of $\Ha^*_{op}$ on $\Ha^\circledast$.
By duality we obtain the following corollary.
\begin{corollary}
The mapping $(\id\tp \tilde \delta_\Ha^l)(\Sc):=\Ru''_{12}\Sc_{13}\Ru_{12}$
defines a structure of left $\Ha$-comodule algebra on $\Ha^\circledast$.
\end{corollary}
Now we can give  one more interpretation of the universal S-matrix.
\begin{thm}
\label{th:intertwiner}
The universal S-matrix $\Sc\in \Ha\tp \Ha^\circledast$ intertwines
the coactions $\tilde\delta_\Ha^l$ and $\delta^l_{\Ha^{op}}$:
\be
\label{l-op-l}
\delta^l_{\Ha^{op}}(a)\Sc=
\Sc\tilde \delta_\Ha^l(a),
\quad
\forall a\in \Ha^\circledast.
\ee
\end{thm}
\begin{proof}
Follows from the reflection equation.
\end{proof}
\end{parag}

\begin{parag}{\bf Bethe subalgebras in $\Ha^*$.}
The remainder of this section is devoted to the so called Bethe subalgebras  in RTT and RE algebras.
Those are maximal commutative subalgebras, which are important in the representation theory
and applications to mathematical physics.
In practice, they are obtained via the construction to be reviewed below,
with the focus on the RE case. That is  a situation when the RE and its conjugate
come along.

We start from the Bethe subalgebras in the dual bialgebra $\Ha^*$ to a
quasitriangular bialgebra $\Ha$.
\begin{definition}
An element  $a\in \Ha^*$ is called invariant if $h\tr a= a\tl h$ for all $h\in \Ha$.
\end{definition}

Denote by $\I\subset \Ha^*$
the subspace of invariant elements.
Clearly $\I$ is a subalgebra, and equation (\ref{RTT}) implies that $\I$ is commutative.
\begin{propn}
\label{I}
Let $\chi\in \Ha$ be a character of the bialgebra $\Ha^*$.
Then the subspace $\I_\chi:=\chi\tr \I=\I\tl\chi$ is a commutative subalgebra
in $\Ha^*$.
\end{propn}
\begin{proof}
An element of $\Ha$ defines a character of $\Ha^*$ if and only if it is group-like.
The left and right actions of $\chi$ are algebra endomorphisms of $\Ha^*$, hence
$\I_\chi$ is a homomorphic image of the commutative subalgebra $\I$.
\end{proof}
The algebra $\I_\chi$ is isomorphic to $\I$, if $\chi$
is invertible in $\Ha$ (e.g. if  $\Ha$ is a Hopf algebra).
We call  $\I_\chi$  a Bethe subalgebra in $\Ha^*$ associated with $\chi$.
Note that we are concerned with the algebraic side of the story and do not address
the practically important question whether $\I_\chi$ is maximal commutative. Than is so
in many interesting cases, see e.g. \cite{NO}.

Next we discuss the Bethe subalgebras in the RE algebras.
\end{parag}
\begin{parag}{\bf Bethe subalgebras in $\Ha^\circledast$.}
Bethe subalgebras in the  algebra $\Ha^\circledast$ are associated with characters of $\Ha^\circledast$.
As $\Ha^\circledast$ coincides with $\Ha^*$ as a vector space,
characters that are identified with elements of $\Ha$.

Every element $\chi\in \Ha$ defines two linear endomorphisms of
$\Ha^*$,
$
a\mapsto \chi\tr a,
$ and $a\mapsto a\tl\chi$.
Consider them as linear endomorphisms of $\Ha^\circledast$.
\begin{propn}
Let $\chi$ be a conjugate character of the algebra $\Ha^\circledast$.
Then the subspace $\I_\chi:=\chi\tr \I=\I\tl \chi$ is a commutative subalgebra in $\Ha^\circledast$.
If two conjugate RE characters $\chi_1,\chi_2$ are invertible in $\Ha$, then the algebras $\I_{\chi_1}$ and $\I_{\chi_2}$ are isomorphic.
\end{propn}
\begin{proof}
For any two elements $a,b\in \I$ one has
$$
\tilde m(a\tl \chi\tp b\tl \chi)= m(a\tl \chi^{(1)}\tp b\tl \chi^{(2)})=m(a\tp b)\tl \chi,
$$
due to (\ref{char_d}) and (\ref{prod}).
It follows immediately  that $\I_\chi$ is a subalgebra and it is a
homomorphic image of the Bethe subalgebra $\I\subset \Ha^*$.
Hence it is commutative and isomorphic to the subalgebra of invariants in $\Ha^*$,
provided $\chi$ is invertible (that is so if $\Ha$ is a Hopf algebra).
\end{proof}
The commutative algebra $\I_\chi$ is called the  Bethe subalgebra of the RE algebra
$\Ha^\circledast$ associated with $\chi$, see \cite{NO}.
\begin{remark}
Note that the isomorphism between $\I_{\chi_i}$ cannot be restricted, in general,
from an automorphism of $\Ha^\circledast$, contrary to the case of $\Ha^*$.
\end{remark}

\end{parag}
\section{Reflection equation and quantum double}
\label{Sec:QD}
\begin{parag}{\bf The T-matrix of a quantum double.}
In this subsection we consider a very special situation when $\Ha$ is the quantum double
of its Hopf subalgebra $\Ha_+\subset \Ha$ and $\Ha_+$ is $\tau$-invariant. For details on this  construction
the reader is referred to \cite{D2}.

Put $\Ha_-=\Ha^{*op}_{+}$.
The double $\D\Ha_+$ is defined as a bialgebra
coinciding with $\Ha_+\tp \Ha_-$ as a vector space.
Both tensor factors are included in $\D\Ha_+$ as sub-bialgebras.
As an $\Ha_+$-$\Ha_-$-bimodule, the double is set to be $\Ha_+\tp \Ha_-$;
the commutation relations between $\Ha_+$ and $\Ha_-$ are expressed through the operations in $\Ha_\pm$ by
\be
\langle\xi^{(1)},x^{(1)}\rangle \xi^{(2)}x^{(2)}=x^{(1)}\xi^{(1)}\langle\xi^{(2)},x^{(2)}\rangle
,\quad
\xi\in \Ha_-,
\quad
x\in \Ha_+.
\label{double_cr}
\ee
We stress that $\langle., .\rangle$ is
the canonical pairing between  $\Ha^*_+$ and $\Ha_+$.

Let $\{e_\mu\}\subset \Ha_{+}$ be a linear basis and let $\{e^\mu\}\subset \Ha_+^*$ be its dual.
\begin{propn}
The canonical element $\Tc\in \Ha\tp \Ha^*$ of the double $\Ha=\D\Ha_+$
 splits to the product
$
(\sum_{\mu}e_\mu \tp e^\mu )(\sum_{\nu}e^\nu\tp e_\nu).
$
\end{propn}
\begin{proof}
We use the fact that $\Ha=\D\Ha_+$ is dual to the twist of the tensor product
$\Ha_+^*\tp \Ha_{+op}$ by the cocycle $I=\sum_\mu e_\mu\tp e^\mu\in \Ha_{+op}\tp \Ha_+^*\subset \Ha\tp \Ha$.
Then the canonical element $\Tc\in \Ha\tp\Ha^*$ is equal to
$\sum_{\mu,\nu}e_\mu\cdot e^\nu\tp e^\mu e^\nu$, where $\cdot$ is the multiplication in the ordinary tensor product
 $\Ha_+\tp \Ha_-$.
Expressing it through the multiplication in $\Ha$ we find $\Tc$ equal to
$$
\sum_{\mu,\nu}(\bar I_1 \tr e_\mu \tl I_{1'})(\bar I_2 \tr e^\nu \tl I_{2'})\tp e^\mu e_\nu=
\sum_{\mu,\nu}e_\mu e^\nu\tp I_{1'} e^\mu \bar I_1 I_{2'} e_\nu \bar I_2.
$$
Here we have used the  symbolic presentation $I=I_1\tp I_2$ and $I^{-1}=:\bar I=\bar I_1\tp \bar I_2$.
We use the dash to distinguish the tensor components of $\bar I$ from the components of  $I$.
Notice that the first  component of $I$ commutes with $\Ha_+^*$ while
the second commutes with $\Ha_{+op}$, in the algebra $\Ha^*$.
The above calculation is continued with $\sum_{\mu,\nu}e_\mu e^\nu\tp e^\mu I_1 \bar I_1e_\nu I_2 \bar I_2$,
which immediately implies the assertion.
\end{proof}
Let $\iota_+$ denote the embedding $\Ha_+\to \Ha$ and $\iota^*_+$ denote the
dual conjugate projection $\Ha^*\to \Ha^*_+$. These maps are Hopf algebra homomorphisms.
\begin{corollary}
The projection $\iota^*_+\colon\Ha^*\to \Ha^*_+$ sends the canonical element $\Tc$ of $\Ha$
to the canonical element of $\Ha_+$, that is, $(\id\tp \iota^*_+)(\Tc)\in \Ha_+\tp \Ha^*_+\subset \Ha\tp \Ha^*_+$.
\end{corollary}

\end{parag}
\begin{parag}{\bf Explicit multiplication in $\Ha^\circledast$ for $\Ha$ being a quantum double.}
\label{expl_mult}
We are going to describe the multiplication in $\Ha^\circledast$ explicitly, assuming that $\Ha$ is a quantum double.
We also assume that not only $\Ha_+$ but also  $\Ha_-$ is $\tau$-invariant. We then say that $\Ha$ is a quantum double
of a Hopf algebra with involution. We show in Section \ref{inf_double} that $\tau$ can be extended to $\Ha$ from
$\Ha_+$, and such an extension is unique under certain conditions.

We choose the triangular structure on $\D\Ha_+$ to be $\Ru=(\sum_\mu e^\mu\tp e_\mu)^{-1}=\sum_\mu\gm^{-1}(e^\mu)\tp e_\mu$
\begin{propn}
\label{Prop:mult_double}
The multiplication in the RE algebra $\Ha^\circledast$ is expressed through the
operations in the Hopf algebra $\Ha^*\simeq \Ha_+^*\tp \Ha_{+op}$ by the formula
\be
\label{mult_expl}
\tilde m(\xi\tp x, \eta \tp y)=
\langle\tau(\eta^{(1)}),x^{(2)}\rangle\bigl(\xi \eta^{(2)}\tp x^{(1)}y\bigr),
\ee
where $\langle . , .\rangle$ is the canonical pairing between   $\Ha^*_+$ and $\Ha_+$.
\end{propn}
\begin{proof}
For all  $\eta\in \Ha_-\subset \Ha$ and $y\in \Ha_+\subset \Ha$
we can write for the coregular action on $\xi\tp x\in \Ha^*_+\tp \Ha_{-}^*=\Ha^*$:
\be
\eta\tr (\xi\tp  x)=\xi \tp \eta\tr x
,\quad  (\xi\tp  x)\tl y=\xi\tl y \tp x.
\label{coreg_gouble}
\ee
These formulas follow from the construction of double as dual to the twisted tensor product
and reflect the  $\Ha_+$-$\Ha_-$-bimodule isomorphism  $\Ha\simeq \Ha_+\tp \Ha_-$.

Explicitly the multiplication (\ref{prod}) can be written as
$$
\tilde m(\xi\tp  x, \eta \tp  y)= \sum_\mu m\Bigl((\tau \gm\gm^{-1})(e^\mu)\tr (\xi\tp x),(\eta\tp  y)\tl e_\mu\Bigr)
= \sum_\mu \xi\bigl(\eta\tl e_\mu\bigr) \tp  (\tau(e^\mu)\tr x) y.
$$
This implies the statement.
\end{proof}

We come to the following immediate  consequence of the Proposition \ref{Prop:mult_double}.
\begin{corollary}
\label{H+-inREa}
(a) The algebras $\Ha_+^*$ and  $\Ha_{+op}$ are embedded in $\Ha^\circledast$ as subalgebras.
(b) As an $\Ha_+^*-\Ha_{+op}$-bimodule, the algebra $\Ha^\circledast$ coincides with $\Ha_+^* \tp \Ha_{+op}$.
\end{corollary}
\end{parag}
\begin{parag}{\bf S-matrix for a double.}
Recall that the Hopf dual $\Ha^*$ to $\Ha=\D\Ha_+$
is isomorphic to the ordinary tensor product $\Ha^*_+\tp \Ha_{-}^*$ of
associative algebras (but not coalgebras). The following result is a corollary of Proposition \ref{Prop:mult_double}.
\begin{thm}
\label{S-factorization}
Suppose that $\Ha$ is the double of a Hopf algebra $\Ha_+$ with involution and let $\{e_\mu\}\subset \Ha_+$
and $\{e^\mu\}\subset \Ha_+^*$ be the dual bases.
Then the universal RE matrix $\Sc\in \Ha\tp \Ha^\circledast$ factorizes to the product $(\sum_\mu e_\mu\tp e^\mu)(\sum_\nu e^\nu\tp e_\nu)$.
\end{thm}
\begin{proof}
The matrix $\Sc$ coincides with the canonical element
$\sum_{\mu,\nu}e_\mu e^\nu\tp e^\mu\cdot e_\nu$ of $\Ha$, where $\cdot$ stands for the multiplication in $\Ha^*$.
The multiplications in $\Ha^*$ and $\Ha^\circledast$ are related by the formula
(\ref{prod}). Thus the canonical element reduces to
$$
\Sc=\sum_{\mu,\nu,\kp} e_\mu e^\nu\tp \bigl((\tau\cdot\gm^{-1})(e^\kp)\tr e^\mu \bigr)(e_\nu\tl e_\kp).
$$
But the formula (\ref{coreg_gouble}) implies $\eta\tr \xi=\xi\ve(\eta)$ and $x\tl y=\ve(y)x$
for all $\eta\in \Ha_-\subset \Ha$, $y\in \Ha_+\subset \Ha$,
$\xi\in \Ha_+^*\subset \Ha^*$, and $x\in \Ha_{+op}\subset\Ha^*$.
This proves the assertion.
\end{proof}
The following result is obtained  under the assumption that $\Ha$ is the double of a Hopf algebra with involution.
\begin{propn}
\label{nochar}
The RE dual $\Ha^\circledast$ has no characters unless  $\dim \Ha=1$.
\end{propn}
\begin{proof}
Suppose, on the contrary, that $\chi\in \Ha$ is a character of $\Ha^\circledast$.
Then, by Corollary \ref{H+-inREa}, it defines
characters $\chi_\pm$ of the algebras $\Ha_\mp$. If follows from Theorem  \ref{S-factorization} that  $\chi=\chi_+\chi_-$ in $\Ha$.
But $\chi_\pm$ are group-like elements in $\Ha$ and so is $\chi$.
The identity $\Delta(\chi)=\chi_1 \Ru'_{12}\chi_2$
implies the equation
$\chi\tp \chi=(\chi\tp 1)\Ru'(1\tp \chi)$. However this equation has no solutions unless $\Ru=1\tp 1$,
since $\chi$ is invertible in $\Ha$ by Proposition \ref{char_inverse}.
\end{proof}
Proposition  \ref{nochar} implies  the following immediate corollary.
\begin{corollary}
If $\Ha$ is a quantum double of a Hopf algebra with involution, then
there exists no algebra homomorphism of $\Ha^\circledast$ to any bialgebra or,
more generally, to any algebra admitting a one-dimensional representation.
\end{corollary}
One of the consequences of this result (rather, its infinite dimensional version) is that the twisted Yangian cannot
be obtained as a quotient of $\Ha^\circledast$ for $\Ha$ being the Yangian double.
\end{parag}
\begin{parag}{\bf Involutions in quantum double.}
\label{inf_double}
We complete the present section  with the
study of involutions in the quantum double. We argue that
any skew involution on $\Ha_+$ can be extended to a skew involution
in a double $\D\Ha_+$. Such an extension is unique if normalized
by the condition $(\tau\tp \tau)(\Ru)=\Ru^{-1}$. Here $\Ru$ is chosen
as in Section \ref{expl_mult}.

Let $\gm$ denote the antipode in $\Ha_-$, which is   the restriction of
the antipode in $\Ha$.
\begin{propn}
\label{inv-double}
Let $\al$ be a $-\choose+$-involution of $\Ha_-$ and $\al^*$ denote the dual
$+\choose-$-involution of $\Ha_+$.
Then the map
$\tau\colon x\xi\mapsto(\al\circ \gm)(\xi)\al^*(x)$, where
$x\in \Ha_+$ and $\xi\in\Ha_-$,
defines a $+\choose-$-involution of $\D\Ha_{+}$.
\end{propn}
\begin{proof}
By construction, the map $\tau$ is an $+\choose-$-involution when restricted to each of the sub-bialgebras
$\Ha_\pm$. So we must check that $\tau$ is an anti-algebra automorphism of $\D\Ha_{+}$,
and it suffices to verify its consistency with
the cross-relations (\ref{double_cr}).
Using the antipode $\gm$ we equivalently rewrite the cross-relations (\ref{double_cr}) in the either form:
$$
\begin{array}{ccc}
\langle\xi^{(1)},x^{(1)}\rangle \xi^{(2)}x^{(2)}\langle\gm(\xi^{(3)}),x^{(3)}\rangle&=&x\xi
\\
\langle\gm(\xi^{(1)}),x^{(1)}\rangle x^{(2)}\xi^{(2)}\langle\xi^{(3)},x^{(3)}\rangle&=&\xi x
\end{array}
,\quad
\xi\in \Ha_-,
\quad
x\in \Ha_+.
$$
To check the consistency, apply the involution $\tau$ to the top line; this will give
the bottom line for $\tau(x)$ and $\tau(\xi)$, as required.
The verification is an easy and straightforward exercise making use of Proposition \ref{sign}.
\end{proof}
\begin{remark}
In general,  the involution $\al\circ \gm$
can be extended to  $\Ha$ in many ways.
Proposition \ref{inv-double} gives  a unique extension
satisfying the equality $\langle \tau (x),\tau(\xi) \rangle=\langle x,\gm(\xi) \rangle$
or, equivalently, $(\tau\tp \tau)(\Ru)=\sum_\mu e^\mu\tp e_\mu=\Ru^{-1}$.
\end{remark}
\end{parag}

\section{Reduced reflection equation algebras}
\label{Sec:FRTform}
\begin{parag}{\bf Completion of  tensor product.}
\label{Subsec:comp_ten_pr}
In this section we put the theory of reflection equation
in the framework of the RTT formalism in its dual version. In this way
we construct RE algebras associated not only
with quasitriangular bialgebras but with their sub-bialgebras
satisfying certain conditions.

We start with a preliminary technical material concerning completion of tensor products.
Consider an associative algebra $\A$  and  its left
module $V$. Let $V^*$ be a submodule in the right $\A$-module $\Hom(V,\C)$ satisfying the following conditions
1) $V\subset \Hom(V^*,\C)$  2) $V^*$ has a basis $\{e^i\}_{i\in I}$ (every element of $V^*$ is a finite linear
combination of $e^i$), 3) there exists the dual system of vectors $\{e_i\}_{i\in I}\subset V$
such that  $\langle e_i,e^i\rangle=\delta_i^j$ for all $i,j\in I$.
The first condition means that the natural pairing between $V$ and $V^*$ is non-degenerate.
We call $V^*$ the restricted dual to $V$.
If $V$ has a restricted dual than we define the completed tensor product
$V\widehat\tp W$ as the vector space of all sums $\sum_{i\in I}e_i\tp w^i$, where $w^i\in W$.
It is a natural $\A$-module  isomorphic to a submodule in $\Hom(V^*,W)$.

Similarly we introduce restricted duals to right $\A$-modules and $\A$-bimodules.
In particular, a restricted dual $\A^*$ to $\A$ considered as a natural bimodule
does exist if $\A$ has a basis and is equipped with a trace-like functional
inducing a non-degenerate inner product.
Then $\A^*$ is isomorphic to $\A$ as a bimodule.
Obviously,
if two algebras $\A$ and $\B$ have restricted duals, so does $\A\tp \B$ and
$(\A\tp \B)^*\simeq\A^*\tp \B^*$.

We also need yet another completion procedure that is special for infinite
direct sums of algebras.
Suppose we are given an infinite family of associative unital algebras
$\{\A^{(k)}\}_{k\in I}$. Consider the infinite direct sum
$\A:=\oplus_{k\in \I}\A^{(k)}$. The elements of $\A$ are infinite sums
$a=\sum_{k\in \I} a^{(k)}$, where $a^{(k)}\in \A^{(k)}$.
The multiplications in $\A^{(k)}$ amount to the  multiplication in $\A$,
with  $a^{(k)}a^{(m)}=0$ if $k\not =m$. The unit in $\A$ is the sum
of units of all summands.
Given two families $\{\A^{(k)}\}_{k\in I}$ and $\{\B^{(k)}\}_{k\in J}$ of algebras
we define their completed tensor product  $\A\bar\tp \B$ to be the algebra
$\oplus_{k,m\in I\times J} \A^{(k)}\tp \A^{(m)}$.
\end{parag}

\begin{parag}{\bf Quasitriangular bialgebra $\Sigma_R\Mg$.}
\label{Subsec:Sigma}
With any solution to the Yang-Baxter equation one can associate a
bialgebra, following  \cite{FRT}.
It is convenient to reformulate that construction  in the dual setting.
Fix an associative algebra $\Mg$ with unit and suppose that $\Mg$ has a restricted dual.
We assume that $\Mg$ is integral over its center $\z(\Mg)$,
that is, the equation $\zt h=0$ with $\zt\in \z(\Mg)$ and $h\in \Mg$ implies either $\zt=0$ or $h=0$.

Put $\Sigma\Mg:=\oplus_{k=0}^\infty \Mg^{\ot k}$ to be the direct sum
of algebras consisting of {\em infinite} formal sums $\sum_{k=0}^\infty a^{(k)}$,
where $a^{(k)}$ are elements from $\Mg^{\ot k}$.
Here the zeroth summand  is set to be the ground field.

The algebra $\Sigma\Mg$ is  a bialgebra. The counit $\ve$
vanishes on $\Mg^{\tp k}$ for all $k>0$ and is identical for $k=0$. The comultiplication is defined as follows. Consider
the "tautological" isomorphism $\pi_{m,n}\colon \Mg^{\tp (m+n)}\to \Mg^{\tp m}\tp\Mg^{\tp n}$, where $m$ and $n$
are non-negative integers. We can think of
$\pi_{m,n}$ as an algebra homomorphism  from $\Mg^{\tp (m+n)}$
to $\Sigma\Mg\bar\tp \Sigma\Mg$. The comultiplication $\Delta\colon\Sigma\Mg\to \Sigma\Mg\bar\tp \Sigma\Mg$
is defined by the formula
\be
\label{coprod}
\Delta := \oplus_{k=0}^\infty \bigl(\sum_{m+n=k}\pi_{m,n}\bigr).
\ee
The counit acts by projecting
to the zeroth components, that is, $\ve\colon h\mapsto h^{(0)}$ for all $h\in \Sigma\Mg$.

Let $R\in \Mg\tp \Mg$ be an {\em invertible} element  fulfilling
the Yang-Baxter equation.
Select in  $\Sigma\Mg$ the subalgebra  $\Sigma_R\Mg$ by imposing  the condition
\be
\label{SigmaM}
\Sigma_R\Mg=\{\sum_{k=0}^\infty h^{(k)}\in \Sigma\Mg\>|\> h^{(k)}\in \Mg^{\tp k},\>R_{i,i+1}h^{(k)}=\si_{i,i+1}(h^{(k)})R_{i,i+1}, \> \forall i=1,\ldots, k-1\}.
\ee
Here $R_{i,i+1}$ is $R$ embedded in $\Mg^{\tp k}$ on the $i$- and $i+1$-th sites, and
$\si_{i,i+1}$ is the flip $\Mg^{\tp k}\to \Mg^{\tp k}$ permuting the  $i$- and $i+1$-th sites.
We denote by $\Mg_R^{(k)}$ the intersection $\Sigma_R\Mg\cap \Mg^{\tp k}$, so that
$\Sigma_R\Mg=\oplus_{k=0}^{\infty} \Mg_R^{(k)}$.

Let $1^{(k)}$ denote the unit of the algebra $\Mg^{\tp k}$.
For all pairs $(k,m)$ of non-negative integers define the elements  $\Ru^{(k),(m)}\in \Mg^{\tp k}\tp \Mg^{\tp m}$
by setting   $\Ru^{(k),(m)}=1^{(k)}\tp 1^{(m)}$ if $km=0$
and  $\Ru^{(k),(m)}:=\prod_{i=1}^k\prod_{j=m}^1 R_{i,j}$ otherwise. The products are ordered  from left to right
 by indicating the initial
and final values of indices.
\begin{propn}
\label{bialgebra}
$\Sigma_R\Mg$ is a sub-bialgebra in  $\Sigma_R\Mg$.
It is quasitriangular, with the
universal R-matrix $\Ru:=\sum_{k,m=0}^\infty \Ru^{(k),(m)}\in \Sigma_R\Mg\bar\tp \Sigma_R\Mg$.
\end{propn}
\begin{proof}
Straightforward.
\end{proof}
\end{parag}
Note that any algebra homomorphism $\rho\colon \Ha\to \Mg$ of a bialgebra $\Ha$ uniquely extends to a bialgebra mapping $\Ha\to \Sigma\Mg$.
Furthermore, if $(\rho^{\tp 2}\circ \Delta)(\Ha)\subset \Mg_R^{(2)}$, then $\rho(\Ha)\subset \Sigma_R\Mg$.

\begin{parag}{\bf RTT algebra $\T_R\Mg^*$.}
Suppose $\Mg$ has a restricted dual $\Mg^*$. Define the restricted duals
to tensor powers $\Mg^{\tp k}$   to be $(\Mg^*)^{\tp k}$.
The tensor algebra $\T\Mg^*$ is equipped with a non-degenerate  pairing with $\Sigma\Mg$
extending the natural pairing  between $\Mg^*$ and $\Mg$.
If $\Mg$ is finite dimensional, then $\T\Mg^*$ is a bialgebra, but
that is not the case in general. However, $\T\Mg^*$ is a natural bimodule algebra
over $\Sigma\Mg$, with the actions  being uniquely extended from the two-sided $\Mg$-action on $\Mg^*$.
Thus $\T\Mg^*$ is a restricted dual to  $\Sigma\Mg$.

Let $\{e_i\}$ be a basis in $\Mg$ and $\{e^i\}$ its dual in $\Mg^*$. Denote
by $T:=\sum_i e_i \tp e^i$   the canonical element
in $\Mg\widehat \tp \Mg^*$.
Consider the quotient $\T_R\Mg^*$ of $\T\Mg^*$ by the ideal generated by
the relations $R T_1 T_2=T_2 T_1R$. Note that the products
here are well defined. Indeed, $T_i T_j$ belongs to the space
$(\Mg \tp \Mg)\widehat \tp(\Mg\tp \Mg)^*$,
which a bimodule over $\Mg\tp \Mg$, cf. Section \ref{Subsec:comp_ten_pr}.

By construction, the elements of the vector space $\Mg^*\subset \T_R\Mg^*$ generate $\T_R\Mg^*$.
The latter algebra admits a non-degenerate pairing
with $\Sigma_R\Mg$ and a two-sided action of  $\Sigma_R\Mg$,
as the RTT relations are two-sided invariant.
\end{parag}
\begin{parag}{\bf RE algebra $\S_R\Mg^*$.}
\label{REalg}
Suppose that the algebra $\Mg$ is equipped with an involutive anti-algebra map
$\tau\colon\Mg\to \Mg$. Then $\tau$ naturally extends to  $\Sigma\Mg$
as a $+\choose-$-involution, which we also denote by $\tau$.
Suppose that the involution $\tau$
restricts to the sub-bialgebra  $\Sigma_R\Mg$.
For instance, that is the case when
the matrix $R$ satisfies the condition $(\tau\tp \tau)(R)\sim R_{21}$
or $(\tau\tp \tau)(R)\sim R^{-1}$, where $\sim$ means equality up
to a multiplier from the center of $\Mg^{\tp 2}$.
We also assume that the restricted dual $\Mg^*$ is stable under the dual conjugate
involution $\tau^*$. Then $\tau^*$ extends to  $\T_R\Mg^*$ as an algebra
automorphism.

The quasitriangular bialgebra $\Sigma_R\Mg$ will play the role of $\Ha$.
Applying  the  RE twist $\F\in (\Ha\bar\tp \Ha)^{\bar\tp 2}$ to  $\T_R\Mg^*$ we
construct the RE dual $(\Sigma_R\Mg)^\circledast$, similarly as in Section \ref{Sec:gen}.
The embedding $\Mg^*\subset (\Sigma_R\Mg)^\circledast$ defines
an algebra homomorphism $\T\Mg^* \to (\Sigma_R\Mg)^\circledast$.
This homomorphism is factorized to the composition
$\T\Mg^*\to (\T\Mg^*)^\circledast \to (\Sigma_R\Mg)^\circledast$,
where the middle term is the twist of the tensor algebra.
The rightmost arrow is the projection along the ideal in $(\T\Mg^*)^\circledast$ which coincides
as the vector space with the defining ideal of the RTT algebra $\T_R \Mg^*$ (a general fact from
the twist theory).
The leftmost arrow is an isomorphism and can be easily described.
 The  natural action of $\Ha^{\bar\tp 2}$ on $\Mg^*$ (involving the involution $\tau$)
extends to an   action of $(\Ha^{\bar\tp 2})^{\bar\tp k}$ on $(\Mg^*)^{\tp k}$.
The isomorphism in question is implemented by the action
of the elements $\F^{(k)}\in (\Ha^{\bar\tp 2})^{\bar\tp k}$ on each tensor power $(\Mg^*)^{\tp k}\subset \T \Mg^*$.
The construction of $\F^{(k)}$ for general twisting cocycle is explained in  Section \ref{Ssec:TTS}.

The kernel of the epimorphism $\T\Mg^*\to (\Sigma_R\Mg)^\circledast$ is generated by the relations  (\ref{RE})
on the matrix  $S=\sum_i e_i \tp e^i$ of $\Mg\widehat \tp \Mg^*$ (again a general fact from
the twist theory, see e.g. \cite{M}).
Introducing the notation  $\S_R\Mg^*$ for the quotient of $\T\Mg^*$
by that ideal, we get an isomorphism $\S_R\Mg^*\simeq (\Sigma_R\Mg)^\circledast$.

If $\T_R\Mg^*$ is a bialgebra, then  $\S_R\Mg^*$ is a right  $\T_R\Mg^*$-comodule
algebra, with the coaction
$$
a\mapsto a^{(2)}\tp \tau^*(a^{(1)})a^{(3)},
$$
expressed through the coproduct in $\T_R\Mg^*$.
\end{parag}
\begin{parag}{\bf Reduced RE algebras.}
It is possible in the framework of the RTT formalism, to construct RE algebras
associated with involutive sub-bialgebras in $\Sigma_R\Mg$.
Suppose that $\Bg$ is a subalgebra in $\Mg$ which is stable under $\tau$.
One can define the bialgebra $\Sigma_R\Bg$ by the same conditions as
$\Sigma_R\Mg$, despite $R$ may not belong to $\Bg\tp \Bg$; thus
$\Sigma_R\Bg=\Sigma_R\Mg\cap \Sigma\Bg$.
By construction, $\Sigma_R\Bg$ is a sub-bialgebra in $\Sigma_R\Mg$
and admits the  involution $\tau$.

Denote by $\T_R\Bg^*$ the quotient of $\T_R\Mg^*$ by the annihilator
of $\Sigma_R\Bg$ under the canonical pairing. As the latter is a bialgebra, the annihilator is an ideal; it
is generated by the annihilator $\Bg^\bot\subset \Mg^*$
of  $\Bg$. Suppose that $\T_R\Bg^*$ is a bialgebra (that is always the case if $\Bg$ is finite dimensional).
Then it is  dual to $\Sigma_R\Bg$; the involution $\tau^*$ descends to  an involution $\T_R\Bg^*$
of the signature $-\choose +$.

Denote by $\S_R\Bg^*$ the quotient of the algebra  $\S_R\Mg^*$ over
the ideal generated by $\Bg^\bot$.
Thus both $\T_R\Mg^*$ and $\S_R\Mg^*$ are generated by the same
vector space $\Mg^*/\Bg^\bot$.

The algebra $\S_R\Mg^*$ is a module algebra over \tw{\Sigma_R\Mg}{\Ru_\tau}{\Sigma_R\Mg}
and hence over $\Sigma_R\Mg$ via the diagonal embedding. Therefore it is an algebra over
$\Sigma_R\Bg\subset \Sigma_R\Mg$.
\begin{propn}
\label{SigmaB-inv}
The defining ideal of $\S_R\Bg^*$ in $\S_R\Mg^*$ is $\Sigma_R\Bg$-invariant.
Hence $\S_R\Bg^*$ is a module algebra over $\Sigma_R\Bg$.
\end{propn}
\begin{proof}
By construction,  $\Sigma_R\Bg$ is annihilated by  $\Bg^\bot$ via the pairing $\langle.,.\rangle$.
For all $\bt\in \Bg^\bot$ and all $a,b\in \Bg$ we have $\langle \Delta(a)\vdash \bt,b\rangle =\langle \bt,\tau (a^{(1)})ba^{(2)}\rangle =0$,
that is, $\Bg^\bot$ is  $\Sigma_R\Bg$-invariant. Hence  $\Bg^\bot$ generates
an invariant ideal in $(\Sigma_R\Mg)^\circledast\simeq \S_R\Mg^*$. The action of  $\Sigma_R\Bg$ on $\S_R\Mg^*$
descends to an action on  $\S_R\Bg^*$.
\end{proof}

As $\T_R\Bg^*$ is a bialgebra  by assumption, we can state the dual version of Proposition \ref{SigmaB-inv}; it can be checked
directly using the defining relations.
\begin{propn}$\S_R\Bg^*$ is a right comodule algebra over the bialgebra $\T_R\Bg^*$.
The coaction on the generators  is expressed through the comultiplication in  $\T_R\Bg^*$
 by the formula
$a\mapsto a^{(2)}\tp \tau^*(a^{(1)})a^{(3)}$, $a\in \Mg^*/\Bg^\bot$.
\end{propn}
Here the restricted coaction on $\Mg^*/\Bg^\bot$ is expressed through the comultiplication in $\T_R\Bg^*\supset \Mg^*/\Bg^\bot$.
\end{parag}
\begin{parag}{\bf The universal S-matrix in the FRT formalism.}
Recall that we have fixed two dual bases $\{e_\mu\}\subset \Sigma\Mg$ and $\{e^\mu\}\subset \T\Mg^*$.
The canonical element $\sum_{\mu}e_\mu \tp e^{\mu}\in \Sigma\Mg\widehat\tp \T\Mg^*$ is a universal T-matrix
as it satisfies the equation (\ref{charateristicT}).

Denote by $\theta$ the projection of $\T\Mg^*\to \T_R\Mg^*$ and
define the matrix $\Tc:=\sum_{\mu}e_\mu \tp\> \theta(e^{\mu})$.
It is clear that  $\Tc$ lies, in fact,  in $\Sigma_R\Mg\widehat \tp \T_R\Mg^*$
and is the universal T-matrix of $\T_R\Mg^*$.

The same canonical element  now considered as that
from $\Sigma_R\Mg\widehat\tp (\Sigma_R\Mg)^\circledast$ is a universal S-matrix, as it satisfies the
equation (\ref{characteristic}). Denote by $\Sc\in \Sigma_R\Mg\widehat\tp \S_R\Mg^*$ its image
under the isomorphism $(\Sigma_R\Mg)^\circledast \to \S_R\Mg^*$.
Using the  equation (\ref{characteristic}), one can recover $\Sc$
from its $\Mg$-component $\Sc^{(1)}$.
\begin{propn}
\label{Prop:fused}
The universal S-matrix equals
\be
\label{fused}
\Sc=\sum_{k=0}^\infty \Sc^{(k)},
\quad \mbox{where}\quad
\Sc^{(k)}=\prod_{i=1}^k\bigl(\Sc^{(1)}_i\prod_{j=i+1}^kR'_{ij}\big).
\ee
The products here are taken in the ascending order from left to right.
\end{propn}
\begin{proof}
Let $\pi^{(k)}$ denote the projection homomorphism $\Sigma\Mg\to \Mg^{\tp k}$ and put $\pi=\pi^{(1)}$.
It is easy to check, using the formula for the coproduct (\ref{coprod}), that $\pi^{(k)}(h)=(\pi^{\tp k}\circ\Delta^k)(h)$ for any $h\in \Sigma\Mg$.
This observation along with (\ref{characteristic}) gives the formula for $\Sc^{(k)}=(\pi^{(k)}\tp \id)(\Sc)$.
\end{proof}
Moreover, one can check that
$\Sc^{(k)}=\Sc^{(i)}_1(\Ru')^{(i),(j)}_{12}\Sc^{(j)}_2$,  where  $i+j=k$. This presentation is independent on
the partition of $k$, which fact follows from coassociativity of comultiplication.

The universal S-matrix for the reduced RE algebra is obtained by  further projection of the right factor to $\S_R\Bg^*$.
Now we consider the inverse problem. Suppose that $\A$ is an associative algebra and suppose
that an element $S\in \Bg\widehat \tp \A$  satisfies the reflection equation (\ref{RE}).
\begin{propn}
The element $\Sc:=\sum \Sc^{(k)}$ with $\Sc^{(k)}=\prod_{i=1}^k\bigl( S_i\prod_{j=i+1}^kR'_{ij}\big)$
is a universal S-matrix with coefficients in $\A$.
\end{propn}
\begin{proof}
Let us demonstrate that $\Sc$ satisfies the equation (\ref{characteristic}).
As the vector space $\Mg^*$ generates $\S_R\Mg^*$, the element $S$ defines a homomorphism $\psi\colon\S_R\Mg\to \A$
by $\psi(a)=\langle a,e_i\rangle s^i$, where $S=\sum_i e_i  \tp s^i$ and $a\in \Mg^*$.
Applying this homomorphism to the equation (\ref{characteristic})
and taking into account (\ref{fused}) with the equality $(\id\tp\psi)(\Sc^{(1)})=S$
we prove the statement.
\end{proof}
\end{parag}
\section{Reflection equation and pseudo-quasitriangular bialgebras}
\label{Sec:Pseudo}
In the present section we relax the condition on the universal R-matrix to be invertible,
in order to incorporate the twisted Yangian into the general picture.
\begin{parag}{\bf Pseudo quasitriangular bialgebras.}
\label{PQHA}
We call an element  $a$ of an associative algebra $\A$ quasi-invertible if
there exists an element $\bar a\in \A$ such that $\bar a a=a\bar a=\zt_a$, where
$\zt_a$ is central in $\A$ and not a zero divisor.
The elements $\bar a$ are $\bar \zt_a$ are not unique, as they can be multiplied by
central elements that are not zero divisors, but we assume that they are fixed. Then we can put  $\zt_a=\zt_{\bar a}$ and $\bar{\bar a}=a$.
The element $\bar a$ will be called quasi-inverse to $a$.
Note that if $\A$ is finite dimensional and has unit, then a quasi-invertible element is always invertible.

We call a bialgebra $\Ha$  pseudo-quasitriangular
if there exists a quasi-invertible element $\Ru\in \A\tp \A$  such that
\be
&(\Delta\tp \id)(\Ru)=\Ru_{13}\Ru_{23},\quad (\id \tp \Delta)(\Ru)=\Ru_{13}\Ru_{12},\nn\\
&(\Delta\tp \id)(\bar\Ru)=\bar\Ru_{23}\bar\Ru_{13},\quad (\id \tp \Delta)(\bar\Ru)=\bar\Ru_{12}\bar\Ru_{13},\nn
\ee
and $\Ru\Delta(h)=\Delta_{op}(h)\Ru$ for all $h\in \Ha$.
It follows that the element
$
(\Delta^k\tp \Delta^m)\Ru
$
is quasi-invertible
in $\Ha^{\tp k}\tp \Ha^{\tp m}$
with the quasi-inverse
$
(\Delta^k\tp \Delta^m)(\bar\Ru).
$

Similarly to the pseudo-quasitriangular structure we define pseudo-twist.
Namely, a quasi-invertible element $\F\in \Ha$ is called
pseudo-twist if it satisfies
the identities
\be
&(\Delta\tp \id)(\F)\F_{12}=(\id\tp \Delta)(\F)\F_{23},
\quad
\bar\F_{12}(\Delta\tp \id)(\bar\F)=\bar\F_{23}(\id\tp \Delta)(\bar\F),
\nn\\
&(\ve\tp \id)(\F)=(\id\tp \ve)(\F)=1\tp 1=(\ve\tp \id)(\bar\F)=(\id\tp \ve)(\bar\F).
\ee
Note that, in general, there is no bialgebra structure on $\Ha$ with the comultiplication $\tilde \Delta$ obeying
$$
\F\tilde \Delta (h)=\Delta (h)\F,
$$
contrary to the case of invertible $\F$.

Given a pseudo-twist and an $\Ha$-module algebra $\A$, still one can construct
its twist $\tilde \A$.
There are two principal differences with the invertible case.
First of all, the  comultiplication $\tilde \Delta$ cannot be defined in general,
so the multiplication in $\tilde \A$ is not compatible with the action
of $\Ha$.
Secondly, if  $\A$ is generated by some $\Ha$-submodule, that is no longer the case for $\tilde \A$.
\end{parag}
\begin{parag}{\bf RE pseudo-twist.}
\label{REa_pseudotwist}
Let $\Ha$ be  a pseudo-quasitriangular bialgebra with involution and $\Ha^*$ its restricted dual.
We assume that the element $\zt_\Ru$ satisfies the identity $(\tau\tp \tau)(\zt_\Ru)=(\zt_\Ru)_{21}$.
Then the central element $\zt'_\Ru=(\tau\tp \id)(\zt_\Ru)$ is symmetric (stable under the flip of tensor factors).
Note that $\zt'_\Ru$ is not a zero divisor.

Let us perform the pseudo-twist of the tensor square $\Ha\tp \Ha$ with the cocycle $\Ru_\tau:=(\tau\tp\tau)(\bar \Ru)$
and construct the corresponding RE dual $\Ha^\circledast$ out of $\Ha^*$.
In the algebra $\Ha^\circledast$ the identity (\ref{REu'}) still holds for any pair of elements $a$ and $b$.
To see that, substitute $\zt_\Ru\vdash(a\tp b)$ into (\ref{RTT}), instead of $a\tp b$.
Then  use
the relation $m\circ \zt_\Ru=m\circ \F\bar\F=\tilde m\circ\bar\F$ between the twisted and non-twisted  multiplications,
centrality of $\zt_\Ru$ and the symmetry of $\zt'_\Ru$.

The characteristic equation for the canonical element $\hat \Sc:=\sum_\mu e_\mu \tp e^\mu\in \Ha\tp \Ha^\circledast$ now reads
\be
\label{pseudo characteristic}
\zt_\Ru'\Delta \hat \Sc=\hat \Sc_1\Ru_{12}' \hat \Sc_2.
\ee
Indeed, multiply (\ref{charateristicT}) by $\zt_\Ru'$ and rewrite the right-hand side
as
$\bigl((\zt'_\Ru)_1\tr\Tc\bigr)\bigl(\Tc\tl(\zt'_\Ru)_2\bigr)$,
using the assumption that $\zt_\Ru$ is central.
Expressing the multiplication through the twisted one and replacing $\Tc$ by $\hat \Sc$ one obtains
the right-hand side of (\ref{pseudo characteristic}).
The identity (\ref{pseudo characteristic}) implies the RE (\ref{uRE}), as
the factor  $\zt'_\Ru=(\zt'_\Ru)_{21}$ can be canceled.
\end{parag}
\begin{parag}{\bf Universal S-matrix.}
Apply the FRT formalism assuming $R\in \Mg\tp \Mg$ to be quasi-invertible.
The universal R-matrix is obtained by the same fusion procedure as in Section \ref{Subsec:Sigma}.
The algebras $\Sigma_R\Mg$ and $\T_R\Mg^*$ are constructed similarly to the invertible case, but now $\Sigma_R\Mg$ is pseudo-quasitriangular.
Consider the RE pseudo-twist of the algebra $\T_R\Mg^*$ with the cocycle $(\tau\tp \tau)(\bar \Ru)$ instead of $\Ru_\tau=(\tau\tp \tau)(\Ru^{-1})$ and denote it $(\Sigma_R\Mg)^\circledast$.
We reserve the notation $\S_R\Mg^*$ for the RE algebra generated by $\Mg^*$ modulo
the relations (\ref{RE}) on the matrix coefficients of $S=\sum_ie_i \tp e^i\in \Mg\widehat{\tp}\Mg^*$.
As argued in the previous subsection, there exists a homomorphism
$\varpi\colon\S_R\Mg^*\mapsto (\Sigma_R\Mg)^\circledast$. However, contrary
to the invertible case, it is not an isomorphism.

The canonical element $\hat \Sc\in \Sigma_R\Mg\widehat \tp (\Sigma_R\Mg)^\circledast$ satisfies
the identity (\ref{pseudo characteristic}), and the formula (\ref{fused}) is generalized to
\be
\label{RE3}
\omega_{\zt'_\Ru} \hat \Sc=
\sum_{k=0}^\infty \prod_{i=1}^k\bigl(\hat \Sc^{(1)}_i\prod_{j=i+1}^kR'_{ij}\big),
\ee
where the element $\omega_{\zt'_\Ru}\in \Sigma_R\Mg$ is computed by the formula
$$
\omega_{\zt'_\Ru}=\sum_{k=0}^\infty \prod_{1\leqslant i<j\leqslant k} (\zt_{R}')_{i,j}.
$$
As $\hat \Sc^{(1)}\in \Mg\widehat{\tp} \varpi(\S_R\Mg^*)$,
the presentation (\ref{RE3}) shows that the element $\Sc:=\omega_{\zt'_\Ru} \hat\Sc$
lies in $\Sigma_R\Mg \widehat{\tp} \varpi(\S_R\Mg^*)$. One can check that that $\Sc$ fulfils (\ref{characteristic}),
thus it is the universal S-matrix with coefficients in  the $\varpi$-image of the algebra $\S_R\Mg^*$.
Unfortunately, the kernel of $\varpi$ is difficult to control. Nevertheless, the right-hand side of (\ref{RE3}) gives a hint how to
define the universal S-matrix with coefficients in  $\S_R\Mg^*$ and for  more general domain of coefficients $\A$.
\begin{propn}
Let $S\in \Mg\widehat{\tp} \A$ be a solution to the reflection equation.
Then the right-hand side of (\ref{RE3}) with $\Sc^{(1)}:=S$
is a universal S-matrix, $\Sc\in \Sigma_R\Mg\widehat{\tp} \A$.
\end{propn}
\begin{proof}
The proof splits into two parts. We show that  $\Sc$ satisfies
(\ref{characteristic}) in $\Sigma\Mg\widehat\tp\A$, and that  $\Sc$ belongs to  $\Sigma_R\Mg\widehat\tp \A$.

First let us check the identity (\ref{characteristic}) in $\Sigma\Mg\widehat{\tp}  \A$.
By construction, the $k$-component of $\Sc$ can be presented as  $\Sc^{(k)}=\Sc^{(1)}_1\Ru^{(1)(k-1)}\Sc^{(k-1)}_{2}$.
By virtue of the formula (\ref{coprod}), condition (\ref{characteristic}) boils down to the fusion identities
$\Sc^{(k)}=\Sc^{(i)}_1\Ru^{(i)(k-i)}\Sc^{(k-i)}_{2}$ for all $k$ and  $i=1,\ldots, k-1$ (plus two obvious cases $i=0,k$). They
 can be proved by the double induction on $k$ and $i$.

Next we must check that $\Sc$ lies, in fact, in $\Sigma_R\Mg\widehat{\tp}  \A\subset \Sigma\Mg\widehat{\tp}  \A$.
That means that the equalities $R_{i,i+1}\Sc^{(k)}=\si_{i,i+1}\Sc^{(k)}R_{i,i+1}$ fulfilled for all $k\geq 2$ and $i=1,\ldots,k-1$.
These equalities are reduced to the case $k=2$ using the mentioned above fusion identities.
\end{proof}

We obtain the following simple but important property of the  universal S-matrices
over the bialgebra $\Sigma_R\Mg$.
\begin{propn}
\label{Prop:homs}
A universal S-matrix with coefficients in an associative algebra $\A$ defines
a homomorphisms $\S_R\Mg^*\to \A$ and vice versa.
\end{propn}
\begin{proof}
Indeed, given a such a homomorphism apply it to the second tensor component of the universal S-matrix of $\Sigma_R\Mg\widehat{\tp}\S_R\Mg^*$
and get an S-matrix with coefficients in $\A$. Conversely, if $\Sc\in\Sigma_R\Mg\widehat{\tp}\A$ is satisfies (\ref{characteristic}), then
its $\Mg$-component $\Sc^{(1)}$ solves the RE in $(\Mg\tp \Mg)\widehat{\tp} \A$. Hence it  defines
a homomorphism  $\S_R\Mg^*\to \A$, as $\S_R\Mg^*$ is generated by $\Mg^*$.
Applying this homomorphism to the S-matrix of $\Sigma_R\Mg\widehat{\tp}\S_R\Mg^*$ we get a universal S-matrix with coefficients in $\A$.
The fusion formula (\ref{characteristic}) guarantees that it coincides with $\Sc$, as it does in the $\Mg$-component.
\end{proof}

Given an involutive subalgebra $\Bg\subset \Mg$,
the RE algebra $\S_R\Bg^*$ is defined as a quotient of $\S_R\Mg^*$ by the ideal generated by $\Bg^\bot$,
the annihilator of $\Bg$ in $\Mg$.
If $\T_R\Bg^*$ is a bialgebra, then $\S_R\Mg^*$ is a right $\T_R\Bg^*$-comodule algebra
under the coaction $S\mapsto \tau^*(T_{13})S_{12}T_{13}$. We call $\S_R\Bg^*$ a {\em reduced} RE algebra,
as in the case of invertible $\Ru$. The universal S-matrix with coefficients
in $\S_R\Bg^*$ is obtained via the projection $\S_R\Mg^*\to \S_R\Bg^*$ from
the universal S-matrix of $\S_R\Mg^*$.
\end{parag}

\section{Application to twisted Yangians}
\label{Sec:TY}
\begin{parag}{\bf The Yangian $\Yg(\g\l_N)$.}
The present section is devoted to applications of the general theory of reflection equation to  twisted Yangians.
We start with the definition of the Yangian $\Yg(\g\l_N)$ within the R-matrix formalism.
Let $\{e_{ij}\}_{i,j=1}^N$ be the standard matrix basis in $\End(\C^N)$ with multiplication
$e_{ij}e_{mn}=\delta_{jm}e_{in}$.
Denote by $P:=\sum_{i,j=1}^N e_{ij}\tp e_{ji}\in\End(\C^{\tp 2N})$
the flip operator.
The Yang R-matrix is a polynomial function of two variables
$u,v$,
\be
\label{R-Yangian}
R(u,v)=(u-v)1_N\tp 1_N-P,
\ee
where $1_N\in \End(\C)$ is the unit matrix.

Denote by $E_{ij}$ the generators of the universal enveloping algebra $\Ug(\g\l_N)$.
They satisfy the relations
$$
E_{ij}E_{km}-E_{km}E_{ij}=\delta_{jk}E_{im}-\delta_{im}E_{kj},
$$
where $i,j,k,m=1,\ldots, N$.
In terms of the operator $E=\sum_{i,j=1}^N e_{ij}\tp E_{ji}$, these relations
can be written as an equation in $\End(\C^N)\tp \End(\C^N)\tp \Ug(\g\l_N)$:
$$
[E_{13},E_{23}]=[E_{23},P_{12}].
$$

The Yangian $\Yg(\g\l_N)$ is an associative algebra generated by the elements $T^{(k)}_{i,j}$,
where $i,j=1,2,\ldots, n$ and $k=1,2,\ldots$. The generators are subject to relations
$$
R(u,v)T_1(u)T_2(v)=T_2(v)T_1(u)R(u,v).
$$
Here $T(u)$ is the matrix $T(u)=1_N+\sum_{k=1}^\infty u^{-k}\sum_{i,j=1}^N e_{ij}\tp T^{(k)}_{ij}
\in \End(\C^N)[[u,u^{-1}]]\tp \Yg(\g\l_N)$.

For a guide in the theory of Yangian $\Yg(\g\l_N)$, the reader is referred to \cite{MNO} or to
the original work \cite{D1}.
\end{parag}
\begin{parag}{\bf The twisted Yangians.}
Denote by $t$  an involutive anti-automorphism of $\End(\C^N)$. This involution factors through
the matrix transposition and conjugation with either symmetric or skew-symmetric matrix. Depending
on its choice, $t$ defines the matrix Lie algebra $\g_N=\{A\in \End(\C^N)|A^t=-A\}$.
Consider the associative
algebra $\Xg(\g_N)$ generated by the
coefficients of the formal series $S_{ij}(u):=\sum_{k=0}^\infty u^{-k} S_{ij}^{(k)}$,
$S_{ij}^{(0)}=\delta_{ij}$,
subject to the relation
\be
R_{12}(u,v)S_1(u)R_{12}'(u,v)S_2(v)=S_2(v)R'_{12}(u,v)S_1(u)R_{12}(u,v).
\quad
\label{re}
\ee
Here
$S(u):=\sum _{i,j=1}^N e_{ij}\tp S_{ij}(u)$ and $R'(u,v)=-(u+v)1_N\tp 1_N-(\tau\tp \id)(P)$. The subscripts
indicate the embeddings
 $\End(\C^N)\to \End(\C^N)\tp \End(\C^N)$.

The algebra $\Xg(\g_N)$ is called the extended twisted Yangian.
It is a right $\Yg(\g\l_N)$-comodule algebra under the coaction
$S(u)\mapsto T^t_{13}(-u)S_{12}(u)\T_{13}(u)$.
The operation $T(u)\mapsto T^t(-u)$ defines an algebra and anti-coalgebra involution on $\Yg(\g\l_N)$,
denoted further by $\tau^*$.

The algebra $\Xg(\g_N)$ admits a homomorphism to the Yangian $\Yg(\g\l_N)$ via the assignment $S(u)\mapsto T^t(-u)T(u)$.
It is a morphism of  $\Yg(\g\l_N)$-comodule algebras, where the coaction on $\Yg(\g\l_N)$ coincides with the comultiplication.
Such a homomorphism is always implemented by a character of $\Xg(\g_N)$, which in our case is specified by the assignment $S_{ij}(u)\mapsto \delta_{ij}$.
The image of this homomorphism is called the twisted Yangian and denoted by $\Yg(\g_N)$.
The involution $\tau^*$ is a coalgebra anti-automorphism of $\Yg(\g\l_N)$.
Hence $\Xg(\g_N)$ becomes as a right  $\Yg^{op}(\g\l_N)$-comodule algebra under the coaction
$S(u)\mapsto T_{12}(u)S_{13}(u)T^t_{12}(-u)$. The assignment $S(u)\to T(u)T^t(-u)$ defines
a map $\Xg(\g_N)\to \Yg^{op}(\g\l_N)$, which is a morphism of right comodules.
The image of this map coincides with the image of $\Yg(\g)$ under the involution $\tau^*$.

\end{parag}
\begin{parag}{\bf Twisted Yangian as a reduced RE algebra.}
\label{Subsec:TY-RREA}
Here we specialize the general constructions of Sections \ref{Sec:gen}, \ref{Sec:FRTform},  and  \ref{Sec:Pseudo} to the case
of twisted Yangians. The role of  $\Ha$ belongs to
the pseudo-quasitriangular bialgebra $\Sigma_R\Mg$ defined by (\ref{SigmaM}), where $\Mg$ is the loop ring,
 $\Mg:=\End(\C^N)[u,u^{-1}]$. The elements $u^{n}e_{ij}$, $n\in \Z$, $i,j=1,\ldots, N$, form a basis in $\Mg$.
The restricted dual $\Mg^*$ is taken to be the linear span of the dual basis $\{t^{-n}_{ji}\}$;
by definition $\langle u^{n}e_{ij},t^{-m}_{kl}\rangle=\delta_{n,m}\delta_{j,k}\delta_{i,l}$.
It is easy to check that $\Mg^*$ is preserved by the coregular two-sided action of $\Mg$.
In fact, $\Mg^*\simeq \Mg$ as a bimodule, via the inner product $\langle Af(u),Bg(u)\rangle=\Tr(AB)\Res_0(uf(u)g(u))$,
where $A,B\in \End(\C^N)$ and $f,g\in \C[u,u^{-1}]$.

Introduce the involutive anti-automorphism of the algebra $\Mg$ by  $\tau\colon f(u)A\mapsto f(-u)A^t$,
for all $f\in \C[u,u^{-1}]$ and $A\in \End(\C^N)$.
The R-matrix (\ref{R-Yangian}) satisfies the Yang-Baxter equation in $\Mg^{\tp 3}$.
It is quasi-invertible, and its quasi-inverse is $\bar R:=(u-v)1_N\tp 1_N +P$ with $\zt_R=(u-v)^2-1$.
Moreover, $(\tau\tp \tau)(R)=R_{21}$, therefore one can construct
the pseudo-quasitriangular bialgebra algebra $\Sigma_{R}\Mg$ with involution, as explained in
Section \ref{Sec:FRTform}.
By induction, one can check that the R-matrix $\Sigma_{R}\Mg$ satisfies the identity $(\tau\tp\tau)(\Ru)=\Ru_{21}$.
From this we conclude that the element the  $\zt_R'$ is symmetric: $(\zt_R')_{21}=\zt_R'$.
Hence we can construct the RE algebra $\S_R\Mg^*$ as in Section \ref{REa_pseudotwist}.

Put $\Bg:=\C1_N+u^{-1}\End(\C^N)[u^{-1}]$, the subalgebra of polynomials currents in $u^{-1}$.
It is obvious that the quotient $\Bg^*=\Mg^*/\Bg^\bot$ is spanned by $\{t^{n}_{ji}\}$, $n\in \N$, $i,j=1,\ldots, N$,
with $t^{0}_{ji}=\delta_{ij}$.
Now construct $\T_{R}\Bg^*$ and the reduced RE algebra $\S_{R}\Bg^*$ as as Section \ref{REa_pseudotwist}. The first algebra coincides with the Yangian $\Yg(\g\l_N)$.
The second algebra coincides with the extended twisted Yangian $\Xg(\g_N)$.

Consider the matrix $\breve R(u,v)=1_N\tp 1_N-\frac{P}{u-v}$ as an element of the completed tensor square of
$\Mg$ via the Laurent expansion
\be
\label{Laurent}
\breve R(u,v)=1_N\tp 1_N-\sum_{k=0}^\infty \frac{v^k}{u^{k-1}}P.
\ee
The assignment $T(u)\mapsto R(u,v)$ or, explicitly, $t^{0}_{ij}\mapsto \delta_{i,j}1_N$,  $t^{m}_{ij}\mapsto v^{m-1}e_{ij}$, for $m\geqslant 1$ and $t^{m}_{ij}\mapsto 0$ for
$m<0$ defines an algebra homomorphism of $\T_{R}{\Mg^*}\to \Mg$. It extends
to a homomorphism $\T_{R}{\Mg^*}\to \Sigma_{R}{\Mg}$, which factors through
a homomorphism $\Yg(\g\l_N)= \T_{R}{\Bg^*}\to \Sigma_{R}{\Mg}$.
It is an anti-coalgebra map.

\end{parag}

\begin{parag}{\bf The double Yangian.}
The double of $\Yg^*(\g\l_N)$ is a natural quasitriangular bialgebra that
contains  $\Yg^*(\g\l_N)$.
Although the twisted Yangians and their S-matrices are related to a different bialgebra, that is, to  $\Ha=\Sigma_{R}\Mg$,
it is interesting to put the double Yangian into  the context.

The quantum double $\D\> \Yg^*(\g\l_N)$ can be described as follows.
It is generated by the elements $L^{\pm,k}_{ij}$, where the indices run over
$i,j=1,2,\ldots, n$ and $k=0,1,\ldots$. The generators are subject to the relations
$$
L^-_1(u)L^-_2(v)(u-v-P_{12})=(u-v-P_{12})L^-_2(v)L^-_1(u),
$$
$$
(u-v-P_{12})L^+_1(u)L^+_2(v)=L^+_2(v)L^+_1(u)(u-v-P_{12}),
$$
$$
L^+_1(v)(v-u-P_{12})L^-_2(u)=L^-_2(u)(v-u-P_{12})L^+_1(v),
$$
were the matrices $L^\pm(u)$ are defined to be
\be
L^+(u)&=&1_N+\sum_{k=0}^\infty u^k\sum_{i,j=1}^N e_{ij}\tp L^{+,k}_{ij},
\nn\\
L^-(u)&=&1_N+\sum_{k=1}^\infty u^{-k}\sum_{i,j=1}^N e_{ij}\tp L^{-,k}_{ij}.
\nn
\ee
The subalgebra generated by the elements $L^{-,k}_{ij}$ can be identified with $\Yg_{op}(\g\l_N)$.
The subalgebra generated by the elements  $L^{+,k}_{ij}$ is identified with $\Yg^{*}(\g\l_N)$.

The coalgebra structure on $\D\> \Yg^*(\g\l_N)$ is introduced on the generators by the formulas
$$
\Delta(L^+)(u)=L^+_{13}(u)L^+_{12}(u),
\quad
\Delta(L^-)(u)=L^-_{12}(u)L^-_{13}(u),
$$
$$
\ve\colon (L^\pm)^k_{ij}\mapsto 0,
$$
This  makes $\D\> \Yg^*(\g\l_N)$ a bialgebra including $\Yg^{*}(\g\l_N)$ and $\Yg_{op}(\g\l_N)$ as sub-bialgebras.
The
pairing between  $\Yg^{*}(\g\l_N)$ and $\Yg(\g\l_N)$ is given on the generators by
\be
\label{pairing}
\langle L^+(u),L^-(v)\rangle=1_N\tp 1_N-P\sum_{k=0}^\infty \frac{v^k}{u^{k+1}}.
\ee
Using the methods of \cite{N}, one can show that this pairing is non-degenerate.
\end{parag}
\begin{parag}{\bf Evaluation homomorphism of the double Yangian.}
It is known that $\Yg_{op}(\g\l_N)$ and hence $\D\> \Yg^*(\g\l_N)$ contains $\Ug(\g\l_N)$ as a sub-bialgebra.
Since $\Ug(\g\l_N)$ is a Hopf algebra, there can be  defined the adjoint action of $\Ug(\g\l_N)$ on $\D\> \Yg^*(\g\l_N)$.
This action makes $\D\> \Yg^*(\g\l_N)$ a $\Ug(\g\l_N)$-algebra.

Denote by $\bar L^-$ the matrix with the elements $\bar L^-_{ij}(u):=\gm^{-1}\bigl(L_{ij}^-(u)\bigr)$, where $\gm$ is the antipode of $\Yg_{op}(\g\l_N)$.
The matrix $\bar L^-$ is the inverse to $L^-$ and can be identified with the generators matrix $T$ of the Yangian $\Yg(\g\l_N)$.
\begin{propn}
\label{module_prop}
a) The coproduct $\Delta\colon \D\> \Yg^*(\g\l_N)\to \D\> \Yg^*(\g\l_N)\tp \D\> \Yg^*(\g\l_N)$ is $\Ug(\g\l_N)$-equivariant.\\
b) The assignment
\be
\label{loop}
\hat\rho \colon \bar L^{-,k}_{ij}\mapsto - z^{k-1}E_{ji},
\quad
\hat\rho\colon L^{+,k}_{ij}\mapsto \frac{1}{z^{k+1}}E_{ji}
\ee
defines a $\Ug(\g\l_N)$-equivariant algebra  homomorphism $\D\> \Yg^*(\g\l_N)$ to the
algebra $\Ug(\g\l_N)[z,z^{-1}]$.\\
c) Under the homomorphism $\rho$ the universal R-matrix $\check\Ru$ of $\D\> \Yg^*(\g\l_N)$ goes over to:
$$
(\rho\tp \id)(\check\Ru)=L^+(z),
\quad
(\rho\tp \id)(\check\Ru_{21})=L^-(z),
$$
$$
(\rho\tp \rho)(\check\Ru)=1_N\tp 1_N-P\sum_{k=0}^\infty \frac{w^k}{z^{k+1}}.
$$

\end{propn}
\begin{proof}
Co-commutativity of $\Ug(\g\l_N)$ implies a).
While checking b), it is convenient to work with the inverse $\bar L^{-}(v)$ of the matrix $L^{-}(v)$; the
defining relations in terms of $\bar L^{-}(v)$ can be readily written down.
Then the matrix form of the assignment (\ref{loop}) is
\be
\bar L^-(u)\mapsto 1_N\tp 1_N-\sum_{k=0}^\infty\frac{z^k}{u^{k+1}}E=1_N\tp 1_N-\frac{E}{u-z},
\nn\\
L^+(u)\mapsto 1_N\tp 1_N+\sum_{k=0}^\infty\frac{u^k}{z^{k+1}}E=1_N\tp 1_N-\frac{E}{u-z}.
\ee
Now b) can be verified using the relation (\ref{loop}).

The statement c) follows from non-degeneracy of the pairing (\ref{pairing}) and the
construction of the double, \cite{D2}.
\end{proof}

The mapping (\ref{loop}) can be considered as a family of homomorphisms $\hat \rho_z\colon \D\> \Yg^*(\g\l_N)\to \End(\C^N)[z,z^{-1}]$,
$z\not=0$.
When restricted to the subalgebra $\Yg(\g\l_N)$, this family extends to $z=0$.
Then $\hat \rho_z$ can be presented as a composition of the
shift automorphism $T(u)\mapsto T(u-z)$ and $\hat \rho_0$.
However, that is not the case in what concerns  the entire algebra  $\D\> \Yg^*(\g\l_N)$.

The basic representation of  $\g\l_N$ induces an algebra homomorphism
$\Ug(\g\l_N)[z,z^{-1}]\to \End(\C^N)[z,z^{-1}]$.
Taking composition of $\hat \rho$ with this homomorphism  we
obtain a representation
$\rho\colon \D\> \Yg^*(\g\l_N)\to \End(\C^N)[z,z^{-1}]$.
As follows from Proposition \ref{bialgebra} and the bottom formula of Proposition \ref{module_prop} c), this representation extends to a bialgebra
homomorphism $\D\> \Yg^*(\g\l_N)\to \Sigma_{R}\Mg$,
with $\Mg:=\End(\C^N)[z,z^{-1}]$.

\end{parag}
\begin{parag}{\bf Universal S-matrix as an intertwiner.}
In this section we use the following modification of $\widehat {\tp }$  that is  special for the case $\Mg=\End(\C^N)[u,u^{-1}]$.
The reason for that is to make the completed tensor product of algebras  an algebra.
Denote by  $\{e_\al\}$ the standard basis in $\Mg^{\tp k}$
(expressed through the powers of the spectral parameters $u_1,\ldots, u_k$).
Given an associative algebra $\A$ we put $\Mg^{\tp k}\widehat \tp \A$ to be the vector space
of all sums $\sum_{\al}e_\al\tp a^\al$, where $a^\al\in \A$ and the powers of the spectral parameters are bounded from above.
The vector space  $\Mg^{\tp k}\widehat \tp \A$ is equipped with the natural structure of associative algebra.
We put
$\Sigma \Mg\>\widehat \tp \A:=\oplus_{k=0}^\infty \Mg^{\tp k}\widehat \tp \A$, the infinite direct sum of algebras.

As we argued in Section \ref{Subsec:Sintertwiner},
the universal S-matrix can be viewed as an intertwiner between representations of the RE algebra.
We are going to prove an analog of the formula (\ref{l-op-l}) for the extended twisted Yangian.

Define the operators $\breve\Ru, \breve\Ru'\in \Sigma_R \Mg\widehat \tp \Sigma_R \Mg$ by setting
\be
\label{R-fusion}
\breve\Ru:=\sum_{k,m=0}^\infty \breve \Ru^{(k),(m)},&
\quad
\breve\Ru^{(k),(m)}:=
\left\{
\begin{array}{ccc}
1^{(k)}\tp 1^{(m)}, &km=0,\\
\prod_{i=1}^k\prod_{j=m}^1 \breve R_{i,j}(u_i,v_j)
, &km\not =0.
\end{array}
\right.\\
\label{R'-fusion}
\breve\Ru':=\sum_{k,m=0}^\infty (\breve \Ru')^{(k),(m)},&
\quad
(\breve \Ru')^{(k),(m)}:=
\left\{
\begin{array}{ccc}
1^{(k)}\tp 1^{(m)}, &km=0,\\
\prod_{i=1}^k\prod_{j=1}^m \breve R'_{i,j}(u_i,v_j)
, &km\not =0.
\end{array}
\right.
\ee
Here the matrix $\breve R_{ij}(u,v)$ is defined by (\ref{Laurent}), and $\breve R'_{ij}(u,v)$
is obtained from $\breve R_{ij}(u,v)$ by the involution $\tau$ applied to the first component.
The matrices $\Ru$ and $\breve \Ru$ have different status in the theory.
While $\Ru$ intertwines the comultiplication of $\Ha$ with its co-opposite,
$\breve \Ru$ implements a homomorphism from $\Yg(\g\l_N)$ to $\Ha$. An analogous interpretation
can be given to $\Ru'$ and $\breve\Ru'$.

Fix a symmetric or skew symmetric numeric $N\times N$-matrix $X$,
that is, obeying the condition $X^t=\pm X$.
It is known that $X$ satisfies
the RE (\ref{re}). In other words, it defines a character via the assignment $S(u)\mapsto X$.
This character can be identified with an element $\chi\in \Sigma_R\Mg$
by the formula
$$
\chi=
\sum_{k=0}^\infty \prod_{i=1}^k\bigl(X_i\prod_{j=i+1}^kR'_{ij}\big).
$$
Consider the following two
algebra maps
\be
\begin{array}{ll}
\varrho^l_\chi\colon \Xg(\g_N)&\stackrel{\chi\tp \delta^r}{\longrightarrow}
\Yg(\g\l_N)\stackrel{\breve \Ru}{\longrightarrow}\Sigma_R\Mg,
\\
\varrho^r_\chi\colon \Xg(\g_N)&\stackrel{\chi\tp \delta^r}{\longrightarrow}
\Yg(\g\l_N)\stackrel{\tau^*}{\longrightarrow}\Yg(\g\l_N)\stackrel{\breve \Ru}{\longrightarrow}
\Sigma_R\Mg
.
\end{array}
\label{varrhos}
\ee
They are, respectively,  homomorphisms
of left and right $\Sigma_R\Mg^*$-comodule algebras.
In terms of the matrices $\Sc\in \Sigma_R\Mg\widehat{\tp} \Xg(\g_N)$ and $\Tc\in \Sigma_R\Mg\widehat{\tp} \Yg(\g\l_N)$
the coaction $\delta^r$ is presented as
$$
\id\tp \delta^r\colon \Sc\mapsto\tau^*(\Tc_{13})\Sc_{12}\Tc_{13}\in \Sigma_R\Mg^-\widehat\tp \bigl(\Xg(\g_N)\tp \Yg(\g\l_N)\bigr),
$$
where $\tau^*(\Tc):=(\id\tp \tau^*)(\Tc)=(\tau\tp \id )(\Tc)$.
Then the maps $\varrho^l_\chi$, $\varrho^r_\chi$ can be written as
as
\be
\id\tp  \varrho^l\colon\Sc\mapsto \breve \Ru'_{12}\chi_1\breve \Ru_{12},
\quad
\id\tp  \varrho^r\colon\Sc\mapsto \breve \Ru_{12}\chi_1\breve \Ru'_{12},
\ee
see Proposition \ref{Prop:homs}.
\begin{propn}
\label{delta}
For all $a\in \Xg(\g\l_N)$ one has
$
\varrho^r_\chi(a)\chi=\chi\varrho^l_\chi(a).
$
\end{propn}
\begin{proof}
There exist central elements
$c,c'$ in  $\Sigma \Mg\widehat \tp \Sigma \Mg$ such that $\Ru=c\breve \Ru$ and $\Ru'=c'\breve \Ru'$.
Being a character of $\Xg(\g_N)$, the element  $\chi$ obeys the reflection equation
$\Ru_{12}\chi_1\Ru_{12}'\chi_2=\chi_2\Ru_{12}'\chi_1\Ru_{12}$,
as an element of $\Sigma_R\Mg$. The proof follows from this equation
via division by $cc'$.
\end{proof}
\end{parag}
\begin{parag}{\bf Concluding remark.}
Many  constructions related to the notion of universal S-matrix, such as the fusion
procedure, Sklyanin determinant, intertwining properties of RE matrices, have been around for quite a time
in representation theory and applications to integrable models, see e.g. \cite{S,MNO,DMS}.
In the present paper we have tried to convert those facts into the general theory of  universal
S-matrix.
The twisted Yangians of Olshanski have served  for us as a working example. While confining ourself to this case, we should state that
the general theory is applicable to the coideal subagebras
of  \cite{No}, and to the q-twisted Yangians of \cite{MRS}.
\end{parag}

\end{document}